\documentclass[12pt]{amsart}
\usepackage[all]{xypic}
\usepackage{amssymb}
\usepackage{amsfonts}
\usepackage{amscd}
\usepackage[english]{babel}
\newtheorem{theorem}{Theorem}[subsection]
\newtheorem{lem}[theorem]{Lemma}
\newtheorem{cor}[theorem]{Corollary}
\newtheorem{prop}[theorem]{Proposition}

\newtheorem{MT}[theorem]{Main Theorem}
\newtheorem{MTS}[theorem]{Motivic Thom-Sebastiani Theorem}
\newtheorem{df-prop}[theorem]{Definition-Proposition}

\theoremstyle{definition}
\newtheorem{definition}[theorem]{Definition}

\theoremstyle{remark}
\newtheorem{remark}[theorem]{Remark}

\theoremstyle{plain}

\theoremstyle{definition}

\theoremstyle{remark}

\theoremstyle{plain}

\def\boxit#1#2{\setbox1=\hbox{\kern#1{#2}\kern#1}%
\dimen1=\ht1 \advance\dimen1 by #1
\dimen2=\dp1 \advance\dimen2 by #1
\setbox1=\hbox{\vrule height\dimen1 depth\dimen2\box1\vrule}%
\setbox1=\vbox{\hrule\box1\hrule}%
\advance\dimen1 by .4pt \ht1=\dimen1
\advance\dimen2 by .4pt \dp1=\dimen2 \box1\relax}

\def\ie{{\emph i.e. \/}}

\renewcommand{\theequation}{\thesubsection.\arabic{equation}}

\def\longhookrightarrow{\mathrel\lhook\joinrel\longrightarrow}

\def\AA{{\mathbf A}}
\def\BB{{\mathbf B}}
\def\CC{{\mathbf C}}

\def\FF{{\mathbf F}}
\def\GG{{\mathbf G}}

\def\LL{{\mathbf L}}

\def\NN{{\mathbf N}}

\def\PP{{\mathbf P}}
\def\QQ{{\mathbf Q}}

\def\ZZ{{\mathbf Z}}

\def\cD{{\mathcal D}}

\def\cF{{\mathcal F}}

\def\cI{{\mathcal I}}

\def\cL{{\mathcal L}}
\def\cM{{\mathcal M}}

\def\cS{{\mathcal S}}

\def\cV{{\mathcal V}}

\mathchardef\alphag="7C0B
\mathchardef\betag="7C0C
\mathchardef\gammag="7C0D
\mathchardef\deltag="7C0E
\mathchardef\varepsilong="7C22
\mathchardef\varphig="7C27
\mathchardef\psig="7C20
\mathchardef\zetag="7C10
\mathchardef\epsilong="7C0F
\mathchardef\rhog="7C1A
\mathchardef\taug="7C1C
\mathchardef\upsilong="7C1D
\mathchardef\iotag="7C13
\mathchardef\thetag="7C12
\mathchardef\pig="7C19
\mathchardef\sigmag="7C1B
\mathchardef\etag="7C11
\mathchardef\omegag="7C21
\mathchardef\kappag="7C14
\mathchardef\lambdag="7C15
\mathchardef\mug="7C16
\mathchardef\xig="7C18
\mathchardef\chig="7C1F
\mathchardef\nug="7C17
\mathchardef\varthetag="7C23
\mathchardef\varpig="7C24
\mathchardef\varrhog="7C25
\mathchardef\varsigmag="7C26
\mathchardef\Omegag="7C0A
\mathchardef\Thetag="7C02
\mathchardef\Sigmag="7C06
\mathchardef\Deltag="7C01
\mathchardef\Phig="7C08
\mathchardef\Gammag="7C00
\mathchardef\Psig="7C09
\mathchardef\Lambdag="7C03
\mathchardef\Xig="7C04
\mathchardef\Pig="7C05
\mathchardef\Upsilong="7C07

\DeclareMathOperator*{\Spec}{Spec}
\begin{document}

\title[Motivic exponential integrals]{Motivic exponential integrals and a motivic Thom-Sebastiani Theorem}

\author{Jan Denef}
\address{University of Leuven, Department of Mathematics,
Celestijnenlaan 200B, 3001 Leuven, Belgium}
\email{Jan.Denef@wis.kuleuven.ac.be}
%\thanks{}

\author{Fran\c cois Loeser}

\address{Centre de Math\'ematiques,
Ecole Polytechnique,
F-91128 Palaiseau
(UMR 7640 du CNRS), {\rm and}
Institut de Math\'{e}matiques,
Universit\'{e} P. et M. Curie, Case 82,
4 place Jussieu,
F-75252 Paris Cedex 05
(UMR 9994 du CNRS)}
\email{loeser@math.polytechnique.fr}
%    General info
%\subjclass{Primary 54C40, 14E20; Secondary 46E25, 20C20}

\date{February 26, 1998}

%\dedicatory{}

%\keywords{Differential geometry, algebraic geometry}

%\begin{abstract}
%This paper is a sample prepared to illustrate the use of the American
%Mathematical Society's \LaTeX{} document class \texttt{amsproc} and
%publication-specific variants of that class.
%\end{abstract}

\maketitle
%\tableofcontents

\renewcommand{\theequation}{\thesection.\arabic{equation}}
\section{Introduction}\subsection{}Let $f$ and $f'$ 
be germs of analytic functions on smooth complex
analytic varieties $X$ and $X'$  and consider the function
$f \oplus f'$ on $X \times X'$ given by
$f \oplus f'  (x, x') = f (x) + f' (x')$. The Thom-Sebastiani
Theorem classically states that the monodromy of $f \oplus f'$ on the
nearby cycles is
isomorphic to the
product of the monodromy of $f$ and the monodromy of $f'$
(in the original form of the Theorem 
\cite{T-S}
the functions
were assumed to have isolated
singularities).
It is now a common idea that the Thom-Sebastiani
Theorem is best understood by using Fourier transformation
and exponential integrals
because of the formula
\begin{equation}\label{1.1}
\int \exp (t (f \oplus f'))
=
\int \exp (t f) \cdot
\int \exp (t f').
\end{equation}
Indeed, by using asymptotic expansions of such integrals
for $t \rightarrow \infty$, A. Varchenko proved a Thom-Sebastiani
Theorem for the Hodge spectrum in the isolated
singularity case
\cite{Varchenko} (see also \cite{Scherk-Steenbrink}) and the general
case has been announced by
M. Saito 
\cite{St2}, \cite{Sa6}.

\bigskip

The aim of the present paper is to give a motivic meaning to equation
(\ref{1.1}) and to deduce
a motivic Thom-Sebastiani
Theorem. To explain our approach, we will begin by reviewing some
known results on $p$-adic exponential integrals.

\subsection{}Let $K$ be a finite extension of $\QQ_{p}$. Let us denote by $R$ the
valuation ring
of $K$, by $P$ the maximal ideal of $R$, and by $k$ the residue field of
$K$. The cardinality of $k$ will be denoted by
$q$, so $k \simeq \FF_{q}$. For
$z$ in $K$, ${\rm ord} \, z \in \ZZ \cup \{+ \infty \}$
denotes the valuation  of $z$, $\vert z \vert  = q^{- {\rm ord} \, z}$,
and ${\rm ac} (z) = z \pi^{- {\rm ord} \, z}$, where $\pi$ is a fixed
uniformizing parameter for $R$.

Let $f \in R[x_{1}, \ldots, x_{m}]$ be a non constant polynomial.
Let $\Phi : R^{m} \rightarrow \CC$ be a locally constant function with
compact support. Let $\alpha$ be
a character of $R^{\times}$, \ie a
morphism
$R^{\times} \rightarrow \CC^{\times}$ with finite image.
For $i$ in $\NN$, set
$$
Z_{\Phi, f, i} (\alpha) =
\int_{\{ x \in R^{m} \, \vert \, 
{\rm ord} f (x) = i\}}
\Phi (x) \alpha ({\rm ac} f (x))|d x|,
$$
where $|d x|$ denotes the Haar measure on $K^{m}$, normalized so that
$R^{m}$
is of measure 1.

We denote by
$\Psi$ the standard additive character on $K$, defined
by
$$z \longmapsto \Psi (z) = \exp (2 i \pi {\rm Tr}_{K / \QQ_{p}}z).$$
For $i$ in $\NN$, we consider the exponential integral
$$
E_{\Phi, f, i} =
\int_{R^{m}}
\Phi (x)
\Psi  (\pi^{- (i + 1)} f (x)) |d x|.
$$
Let $\alpha$ be a character of $R^{\times}$. The conductor of $\alpha$,
$c ({\alpha})$, is defined as
the smallest $c \geq 1$ such that $\alpha$ is trivial
on $1 + P^{c}$, and one associates to $\alpha$ the Gauss sum
$$
g (\alpha) =  q^{1 - c (\alpha)}
\sum_{v \in (R / P^{c (\alpha)})^{\times}}
\alpha (v) \Psi (v / \pi^{c (\alpha)}).
$$

\begin{prop}\label{p-adique}(See \S\kern .15em 1
of \cite{Denef Bourbaki}.)
For any $i$ in $\NN$,
\begin{equation*}
E_{\Phi, f, i} = \int_{\{ x \in R^{m} \, \vert \, 
{\rm ord} f (x) > i\}}
\Phi (x) |d x|
+ (q - 1)^{-1}
\sum_{\alpha} g (\alpha^{-1}) Z_{\Phi, f, i - c(\alpha) + 1} (\alpha).
\end{equation*}
Here $i - c(\alpha) + 1 \geq 0$. If moreover
the critical locus of
$f$ in ${\rm Supp} \, \Phi$
is contained in $f^{-1} (0)$, then,
for all except a
finite number
of characters $\alpha$, the integrals $Z_{\Phi, f, j} (\alpha)$
are zero for all $j$. 
\end{prop}

\begin{cor}\label{1.2.2}(Using Theorem 3.3 of \cite{Denef Bourbaki}.)
Assume that $\Phi$ is residual, \ie that
${\rm Supp} \, \Phi$ is contained in $R^{m}$ and that $\Phi (x) $
depends only on $x$ modulo $P$,  and that the 
critical locus of
$f$ in ${\rm Supp} \, \Phi$
is contained in $f^{-1} (0)$.
Assume furthermore that the divisor $f = 0$ has good reduction
(in the sense that the conditions in Theorem 3.3 of
\cite{Denef Bourbaki} are satisfied).
Then
\begin{equation*}
E_{\Phi, f, i} = \int_{\{ x \in R^{m} \, \vert \, 
{\rm ord} f (x) > i\}}
\Phi (x) |d x|
+ (q - 1)^{-1}
\sum_{\genfrac{}{}{0pt}{}{\alpha}{c (\alpha) = 1}} g (\alpha^{-1}) Z_{\Phi, f, i} (\alpha).
\end{equation*}
\end{cor} 
\renewcommand{\theequation}{\thesubsection.\arabic{equation}}

So we see that $p$-adic exponential integrals may be expressed
as linear combinations of $p$-adic  integrals
involving multiplicative characters
with Gauss sums as
coefficients.

\medskip

When $k$ is a field of characteristic 0, there is a $k ((t))$-analogue
of
$p$-adic integration, motivic integration,
introduced by M. Kontsevich.
In particular, it is possible by the results of
\cite{Motivic Igusa functions}
and
\cite{Arcs} to define motivic analogues of the 
$p$-adic  integrals $Z_{\Phi, f, i} (\alpha)$
in \ref{1.2.2} as elements of a Grothendieck group
of Chow motives. 
Since in  this analogy the
$k ((t))$-case always has good reduction,
it becomes quite natural to use equality
\ref{1.2.2} as a candidate for the definition of
motivic exponential integrals
and this is indeed what we do in this paper. To achieve this aim,
we enlarge
slightly
our virtual motives by attaching virtual motives to Gauss sums
in a way very similar to Anderson's construction of 
ulterior motives \cite{Anderson}. But now
equation
(\ref{1.1}) is no longer trivial, and the main
result of the paper is that it still holds true
for our motivic exponential integrals
(Theorem \ref{MT}).
We deduce from this result a motivic analogue of the Thom-Sebastiani
Theorem
(Theorem \ref{MTS}). Passing to Hodge realization, this gives a proof of
the
Thom-Sebastiani
Theorem for the Hodge spectrum (Corollary \ref{spectrum}).

\section{Adding ulterior motives to the Grothendieck group}

\subsection{}
We fix a base field $k$, which we assume throughout the paper to be of
characteristic zero,
and we denote by
$\cV_{k}$ the category of smooth and projective $k$-schemes.
For  an object $X$ in $\cV_{k}$ and an integer $d$,
we denote by  $A^{d} (X)$ the Chow group of codimension
$d$ cycles with rational coefficients
modulo rational equivalence.
Objects of the category $\cM_{k}$ of (rational)
$k$-motives
are triples $(X, p, n)$ where $X$ is in $\cV_{k}$,
$p$ is an idempotent (\ie $p^{2} = p$) in the ring of
correspondences ${\rm Corr}^{0} (X, X)$
($= A^{d} (X \times X)$ when $X$ is of pure dimension $d$), and
$n$ is an integer. If $(X, p, n)$
and $(Y, q, m)$ are motives, then
$$
{\rm Hom}_{\cM_{k}} ((X, p, n), (Y, q, m))
=
q \, {\rm Corr}^{m - n} (X, Y) \, p.
$$
Composition of morphisms is given by composition of correspondences.
The category $\cM_{k}$ is  additive, $\QQ$-linear, and pseudo-abelian,
and 
there is a natural tensor product on $\cM_{k}$.
We denote by $h$ the functor $h : \cV_{k}^{\circ} \rightarrow
\cM_{k}$ which sends an object $X$ to $h (X) = (X, {\rm id}, 0)$
and a morphism $f : Y \rightarrow X$ to its graph in
${\rm Corr}^{0} (X, Y)$. We denote by ${\LL}$ the Lefschetz motive
$\LL = ({\rm Spec} \, k, {\rm id}, -1)$. There is a canonical isomorphism
$h (\PP^{1}_{k}) \simeq 1 \oplus \LL$.

When $E$ is a field of characteristic zero,
one defines similarly
the category $\cM_{k, E}$ of $k$-motives with coefficients in $E$, by replacing
the Chow groups
$A^{\cdot}$ by $A^{\cdot} \otimes_{\QQ} E$.
Let $K_{0} (\cM_{k, E})$ be the Grothendieck group of the pseudo-abelian
category
$\cM_{k, E}$. It is also the abelian group associated to the monoid of
motives
with coefficients in $E$. The tensor product on $\cM_{k, E}$ induces 
a natural ring structure on $K_{0} (\cM_{k, E})$.
For $m$ in $\ZZ$, let $F^{m} K_{0} (\cM_{k, E})$
denote 
the subgroup of $K_{0} (\cM_{k, E})$ generated by
$h (S, f, i)$, with $i - \dim S \geq m$. This gives a filtration
of the ring $K_{0} (\cM_{k, E})$
and we denote by 
$\widehat K_{0} (\cM_{k, E})$
the completion of 
$K_{0} (\cM_{k, E})$ with respect to this filtration.

\begin{remark}\label{separation}We expect, but do not know how to prove, that
the filtration $F^{\cdot}$ on 
$K_{0} (\cM_{k, E})$
is separated. This assertion is clearly implied by the
conjectural existence (cf. \cite{Sc} p.185) of additive functors
$h^{\leq j} : \cM_{k, E} \rightarrow \cM_{k, E}$,
$j \in \ZZ$, such that for any $X = h (S, f, i)$
in 
$\cM_{k, E}$, the $h^{\leq j} (X)$ form a filtration of $X$ with
$h^{\leq - 2i -1} (X) = 0$, $h^{\leq 2 \dim S - 2 i} (X) = X$, and
$h^{\leq j} (\LL X) = \LL \, h^{\leq j - 2} (X)$
for all $j$.

In particular, using the existence of weight filtrations, 
one obtains, without using any conjecture, that the kernel
of the canonical morphism
$$ K_{0} (\cM_{k, E}) \longrightarrow \widehat K_{0} (\cM_{k, E})$$
is killed by 
\'{e}tale and Hodge realizations.
\end{remark}

\subsection{}We begin by recalling some material from
\cite{Motivic Igusa functions}.
Let $G$ be a finite abelian
group and let $\hat G$ be its complex character group.
We denote by 
$\cV_{k, G}$ the category of smooth and projective $k$-schemes
with $G$-action.
Let $E$ be a subfield of $\CC$
containing all the roots of unity of order dividing $|G|$. 
For $X$ in $\cV_{k, G}$ and $g$ in $G$, we denote by $[g]$ the
correspondence given by the graph of multiplication by $g$.

For $\alpha$ in $\hat G$,
we consider the idempotent 
$$
f_{\alpha} := \vert G \vert^{-1} \sum_{g \in G}
\alpha^{-1} (g)
[g]
$$
in 
${\rm Corr}^{0} (X, X) \otimes E$, and
we denote by
$h (X, \alpha)$ the motive
$(X, f_{\alpha}, 0)$ in $\cM_{k, E}$.
We will denote by ${\rm Sch}_{k, G}$ the category of  
separated schemes of finite type over $k$ with $G$-action
satisfying the following condition:
the $G$-orbit of any closed point of $X$ is contained in
an affine open subscheme. This condition is clearly satisfied
for $X$ quasiprojective and insures the existence of
$X / G$  
as a scheme. Objects of ${\rm Sch}_{k, G}$ will
be called $G$-schemes and in this paper all 
schemes with $G$-action will be assumed to be $G$-schemes.

The following result, proved in
\cite{Motivic Igusa functions} 
as a consequence of results in \cite{G-N} and \cite{V2}, 
generalizes to the $G$-action case a result which has been proved by
Gillet and Soul\'{e} \cite{G-S} and
Guill\'{e}n and Navarro Aznar \cite{G-N}.

\begin{theorem}Let $k$ be a field
of characteristic 0.
There exists a unique map
$$\chi_{c} : {\rm Ob} {\rm Sch}_{k, G} \times \hat G
\longrightarrow K_{0} (\cM_{k, E})$$
such that
\begin{enumerate}
\item[(1)] If $X$ is smooth and projective
with $G$-action, for any character $\alpha$,
$$\chi_{c} (X, \alpha) = [h (X, \alpha)].$$
\item[(2)] If $Y$ is a closed $G$-stable
subscheme in a scheme
$X$ with $G$-action, for any character $\alpha$,
$$\chi_{c} (X \setminus Y, \alpha) = \chi_{c} (X, \alpha) - \chi_{c} (Y, \alpha).$$
\item[(3)] If $X$ is 
a scheme with $G$-action, $U$ and $V$ are $G$-invariant open
subschemes of $X$, for any character $\alpha$,
$$
\chi_c (U \cup V, \alpha)
=
\chi_c (U, \alpha) +
\chi_c (V, \alpha) -
\chi_c (U \cap V, \alpha).
$$
\end{enumerate}
Furthermore, $\chi_{c}$ is determined by conditions
(1)-(2).
\end{theorem}

\subsection{}In this subsection we gather some elementary statements
we shall need.

\begin{prop}\label{2.3.1}
Let $k$ be a field
of characteristic 0.
\begin{enumerate}
\item[(1)]For any $X$ in ${\rm Ob} {\rm Sch}_{k, G}$,
$$
\chi_{c} (X) = \sum_{\alpha \in \hat G} \chi_{c} (X, \alpha).
$$
\item[(2)]Let $X$ be in ${\rm Ob} {\rm Sch}_{k, G}$. Assume the
$G$-action factors through a morphism of finite abelian
groups
$G \rightarrow H$.
If $\alpha$ is not in the image of $\hat H \rightarrow \hat G$,
then $\chi_{c} (X, \alpha) = 0$.
\item[(3)]Let $X$ and $Y$ be in ${\rm Ob} {\rm Sch}_{k, G}$
and let $G$ act diagonally on $X \times Y$.
Then
$$
\chi_{c} (X \times Y, \alpha) = \sum_{\beta \in \hat G} \chi_{c} (X,
\beta) \, \cdot \,
\chi_{c} (Y, \alpha \beta^{-1}).
$$
\item[(4)]Let $X$ be in ${\rm Ob} {\rm Sch}_{k, H}$ and let
$f : G \rightarrow H$ be a morphism of finite abelian groups.
Then for any $\alpha$ in $\hat G$,
$$
\chi_{c} (X, \alpha)
=
\sum_{\hat f (\beta) = \alpha} \chi_{c} (X, \beta),
$$
with $\hat f : \hat H \rightarrow \hat G$
the dual morphism between character groups.
\end{enumerate}
\end{prop}

\begin{proof} Statements
(1)-(3) are proven in Proposition 1.3.3 of \cite{Motivic Igusa
functions}.
They are all consequence of (4), whose proof is similar~: the statement
is obvious when $X$ is smooth projective, and the general case
follows by additivity of $\chi_{c}$.
\end{proof}

\begin{lem}\label{Kummer}Let $a$ be an integer and let
$\mu_{d} (k)$ act on $\GG_{m, k}$ by
multiplication by $\xi^{a}$, $\xi \in \mu_{d} (k)$.
For any non trivial character $\alpha$ of
$\mu_{d} (k)$, $$\chi_{c} (\GG_{m, k}, \alpha) = 0.$$
\end{lem}

\begin{proof} This is Lemma  1.4.3 of \cite{Motivic Igusa functions}.
\end{proof}

We now
discuss motivic Euler
characteristics of quotients.
The following lemma is well known.

\begin{lem}\label{lem-quot}Let $X$ be a smooth projective
scheme with $G$-action, $H$ a subgroup of $G$, and
$\alpha$ a character of $G / H$. Assume
the quotient $X / H$ is smooth. Then
$h (X / H, \alpha) \simeq h (X, \alpha \circ \varrho)$,
where $\varrho$ is the projection $ G \rightarrow G / H$.
\end{lem}

\begin{proof} This is Lemma  1.5.1 of \cite{Motivic Igusa functions}.
\end{proof}

In general one has the following result, which was conjectured in
\cite{Motivic Igusa functions} and proved in \cite{B-N} by
del Ba{\~ n}o Rollin and Navarro Aznar.

\begin{theorem}\label{assertion}
If $X$ is 
a scheme with $G$-action, $H$ a subgroup of $G$, and
$\alpha$ a character of $G / H$, then
$\chi_{c} (X / H, \alpha) = \chi_{c} (X, \alpha \circ \varrho)$,
where $\varrho$ is the projection $ G \rightarrow G / H$.
\end{theorem}

\begin{remark}\label{prop-quot}Theorem \ref{assertion} 
is a direct consequence of Lemma \ref{lem-quot}
when $X$ is a smooth curve or
when $X$ is smooth and may be embedded in 
a smooth projective
scheme $Y$ with $G$-action such that the quotient $Y / H$ is smooth,
and such that $Y \setminus X$ is the union of finitely many 
smooth closed $G$-stable subvarieties intersecting transversally and
having smooth images in $Y / H$
which intersect transversally.
\end{remark}

\subsection{Jacobi motives}\label{Jacobi}We fix an
integer $d \geq 1$. We denote by $\mu_{d} (k)$
the group of $d$-roots of 1 in $k$ and by $\zeta_{d}$ a fixed primitive
$d$-th root of unity in $\CC$.
We assume from now on that $k$
contains all $d$-roots of unity.

We set
$$ A_{d}    := K_{0} (\cM_{k, \QQ[\zeta_{d}]}).$$

For $n \geq 1$, we consider the affine Fermat
variety $F^{n}_{d}$ defined by
the equation
$x_{1}^{d} + \cdots  + x_{n}^{d} = 1$ in $\AA^{n}_{k}$,
and its closure in  $\PP^{n}_{k}$, which we denote by
$W^{n}_{d}$. Hence 
$W^{n}_{d}$ is defined in $\PP^{n}_{k}$
by
the equation
$- X_{0}^{d} + \cdots  + X_{n}^{d} = 0$,
with $x_{i} = X_{i} / X_{0}$, $i \geq 1$.

The action of $\mu_{d} (k)$ on each coordinate induces a natural action of
the group $\mu_{d} (k)^{n}$ on $F^{n}_{d}$.
Hence, for $\alpha_{1}, \ldots, \alpha_{n}$ characters
of $\mu_{d} (k)$, one defines the Jacobi motive
$J (\alpha_{1},
\ldots, \alpha_{n})$ as the 
element
\begin{equation*}
J (\alpha_{1},
\ldots, \alpha_{n})
:=
\chi_{c} (F^{n}_{d}, (\alpha_{1},
\ldots, \alpha_{n})) 
\end{equation*}
in
$A_{d}$. It is clear that 
$J (\alpha_{1},
\ldots, \alpha_{n})$ is symmetric in the $\alpha_{i}$'s.
We also define
$[\alpha (- 1)] := \chi_{c} (x^{d} = - 1, \alpha)$.
Remark that $[\alpha (- 1)] = 1$, if $k$ contains a $d$-th root of -1.
We will need the following proposition, which is classical in other contexts.

\begin{prop}\label{formulaire}The following relations hold in
$A_{d}$.
\begin{enumerate}
\item[(1)]We have $J (1, 1) = \LL$.
\item[(2)]We have
$J (1, \alpha) = 0$ if $\alpha \not = 1$.
\item[(3)]If $\alpha \not = 1$,
$ J (\alpha, \alpha^{-1}) = - [\alpha (- 1)].$
\item[(4)]We have \begin{equation*}
J (\alpha_{1}, \alpha_{2})
[J (\alpha_{1} \alpha_{2}, \alpha_{3}) - \varepsilon]
=
J (\alpha_{1}, \alpha_{2}, \alpha_{3}) - \delta,
\end{equation*}
with $\varepsilon = \delta = 0$ if $\alpha_{1} \alpha_{2} \not=1$,
$\varepsilon = 1$, $\delta = [\alpha_{1} (- 1)] (\LL - 1)$, if
$\alpha_{1} \alpha_{2} =1$
and $\alpha_{1} \not=1$, and $\varepsilon = 1$, $\delta = \LL$, if $\alpha_{1} = \alpha_{2}=1$.
\end{enumerate}
\end{prop}

\begin{proof}Relation (1) follows directly
from Remark \ref{prop-quot}
because the quotient of the curve
$F^{2}_{d}$ by $\mu_{d} (k) \times \mu_{d} (k)$
is the affine line $\AA^{1}_{k}$.

To prove (2), observe that the quotient of 
the curve
$F^{2}_{d}$ by $\mu_{d} (k) \times \{1\}$ is the Kummer cover
$x_{2}^{d} = 1 - x_{1}$ of the affine line.
Hence, by Remark \ref{prop-quot},
$J (1, \alpha) = \chi_{c} (\AA^{1}_{k}, \alpha)$,
with the natural action of $\mu_{d} (k)$
on
$\AA^{1}_{k}$, and the assertion follows from
Lemma \ref{Kummer} and Proposition \ref{2.3.1}.

Let us now prove (3). Observe that, the character
$\alpha$ being non trivial, we have
$$\chi_{c} (F^{2}_{d} \cap \{x_{2} = 0\}, (\alpha, \alpha^{-1}))=0.$$
Indeed, this follows for instance from Proposition \ref{2.3.1}.
Hence $$J (\alpha, \alpha^{-1}) =
\chi_{c} (F^{2}_{d} \setminus \{x_{2} = 0\}, (\alpha, \alpha^{-1})).$$
Now we may identify
$F^{2}_{d} \setminus \{x_{2} = 0\}$ with the affine curve
$u^{d} - v^{d} = -1 ; v \not= 0$, via the change of variable
$u = x_{1} x_{2}^{-1}$, $v = x_{2}^{-1}.$
Taking in account the 
$\mu_{d} (k) \times \mu_{d} (k)$-action, we get
$$
\chi_{c} (F^{2}_{d} \setminus \{x_{2} = 0\}, (\alpha, \alpha^{-1}))
=
\chi_{c} ((u^{d} - v^{d} = -1 ; v \not= 0), (\alpha, 1)).
$$
By Remark \ref{prop-quot},
\begin{align*}
\chi_{c} ((u^{d} - v^{d} = -1 ; v \not= 0), (\alpha, 1))
&=
\chi_{c} ((u^{d} = v -1 ; v \not= 0), \alpha)\\
& = \chi_{c} (\AA^{1}_{k}, \alpha) - \chi_{c} (u^{d} = - 1, \alpha)\\
& = -[\alpha (- 1)],
\end{align*}
because $\chi_{c} (\AA^{1}_{k}, \alpha) = 0$ by 
Lemma \ref{Kummer} and Proposition \ref{2.3.1}.

To prove (4) we will consider the morphism
$$f : F^{2}_{d} \times (F^{2}_{d} \setminus \{y_{1} = 0\})
\longrightarrow F^{3}_{d} \setminus \{z_{1}^{d} + z^{d}_{2} = 0\}$$
given by
$$
((x_{1}, x_{2}), (y_{1}, y_{2})) \longmapsto (z_{1}
= x_{1} y_{1}, z_{2} = x_{2} y_{1}, z_{3} = y_{2}).
$$
This morphism identifies 
$F^{3}_{d} \setminus \{z_{1}^{d} + z^{d}_{2} = 0\}$
with the quotient of 
$F^{2}_{d} \times (F^{2}_{d} \setminus \{y_{1} = 0\})$
by the kernel $\Gamma$
of the morphism
$\mu_{d} (k)^{4} \rightarrow \mu_{d} (k)^{3}$
given by $(\xi_{1}, \xi_{2}, \xi_{3}, \xi_{4}) \mapsto
(\xi_{1}\xi_{3}, \xi_{2}\xi_{3}, \xi_{4})$.
It follows from
Remark \ref{prop-quot} and Proposition \ref{2.3.1} (3)
that
\begin{multline*}
\chi_{c} (F^{2}_{d}, (\alpha_{1}, \alpha_{2})) \cdot
\chi_{c} (F^{2}_{d} \setminus \{y_{1} = 0\},
(\alpha_{1} \alpha_{2}, \alpha_{3}))\\
=
\chi_{c} (
F^{3}_{d} \setminus \{z_{1}^{d} + z^{d}_{2} = 0\},
(\alpha_{1}, \alpha_{2}, \alpha_{3})).
\end{multline*}
Indeed, $f$ extends to the rational map
$$W^{2}_{d} \times W^{2}_{d}
\longrightarrow W^{3}_{d}$$
given by
\begin{multline*}
([X_{0}, X_{1}, X_{2}], [Y_{0}, Y_{1}, Y_{2}]) \longmapsto\\ [Z_{0} =
X_{0} Y_{0},
Z_{1}
= X_{1} Y_{1}, Z_{2} = X_{2} Y_{1}, Z_{3} = X_{0} Y_{2}].
\end{multline*}
Now let $Z$ and $Z'$ denote respectively
the blow-up of 
$W^{2}_{d} \times W^{2}_{d}$ along $\{X_{0} = 0\} \times \{Y_{1} = 0\}$
and the blow-up of $W^{3}_{d}$ along $\{Z_{0} = Z_{3} = 0\} \cup
\{Z_{1} = Z_{2} = 0\}$. It is classical (see
\cite{Shioda}) that $f$ extends to a morphism
$\bar f : Z \mapsto Z'$, that the actions of
$\mu_{d} (k)^{4}$ and
$\mu_{d} (k)^{3}$ extend to actions on $Z$ and $Z'$ respectively, and
that
$\bar f$ identifies $Z'$ with the quotient of $Z$ by $\Gamma$. (Of course
one could also use
here Theorem \ref{assertion} directly
instead of 
Remark \ref{prop-quot}.)
To complete the proof of (4), we only have to prove
that $\chi_{c} (F^{2}_{d} \cap \{y_{1} = 0\},
(\alpha_{1} \alpha_{2}, \alpha_{3})) = \varepsilon$
and that 
$\chi_{c} (
F^{3}_{d} \cap \{z_{1}^{d} + z^{d}_{2} = 0\},
(\alpha_{1}, \alpha_{2}, \alpha_{3})) = \delta$.
The first equality is clear since, by Proposition
\ref{2.3.1},
$$
\chi_{c} (F^{2}_{d} \cap \{y_{1} = 0\},
(\alpha_{1} \alpha_{2}, \alpha_{3})) =
\chi_{c} (y_{1} =
0, \alpha_{1} \alpha_{2}) \cdot \chi_{c} (y_{2}^{d} =1
, \alpha_{3}).
$$
To prove the second one, we remark that
\begin{multline*}\chi_{c} (
F^{3}_{d} \cap \{z_{1}^{d} + z^{d}_{2} = 0\} \cap
\{z_{1} \, \hbox{or} \, z_{2} \not= 0\},
(\alpha_{1}, \alpha_{2}, \alpha_{3})) =\\
\chi_{c} (u^{d} = -1, \alpha_{1}) \cdot \chi_{c} (\GG_{m, k}, \alpha_{1}
\alpha_{2}) \cdot \chi_{c} (w^{d} = 1, \alpha_{3}).
\end{multline*}
This follows from 
Proposition
\ref{2.3.1}, by using the change of variable
$u = z_{1} z_{2}^{-1}$, $v = z_{2}$ and $w = z_{3}$.
The result is now a consequence of
Lemma \ref{Kummer}, since
$\chi_{c} (
F^{3}_{d} \cap \{z_{1}^{d} + z^{d}_{2} = 0\} \cap
\{z_{1} = z_{2} = 0\},
(\alpha_{1}, \alpha_{2}, \alpha_{3}))$ is equal to 1 or 0, according
whether $\alpha_{1}$ and $\alpha_{2}$ are both trivial or not.
\end{proof}

\subsection{}We assume from
now on that $k$ contains all the roots of unity.
When $d$ divides $d'$ we have a canonical surjective
morphism of groups $\mu_{d'} (k) \rightarrow \mu_{d} (k)$
given by $x \mapsto x^{d' / d}$ which dualizes to a injective morphism
of character groups
$\widehat \mu_{d} (k) \rightarrow \widehat \mu_{d'} (k)$.
We set $\widehat \mu (k) := \varinjlim \widehat
\mu_{d} (k)$. We shall identify $\widehat \mu_{d} (k)$
with the subgroup of elements of order dividing $d$ in $\widehat \mu 
(k)$.

We denote by $F$ the subfield of $\CC$ generated by the roots of unity
and we set
$$
A := K_{0} (\cM_{k, F})$$
and
$$
\widehat A := \widehat K_{0} (\cM_{k, F}).$$
We have a natural ring morphism
$A_{d} \rightarrow A$. When no confusion occurs, 
we will still denote by the same symbol the image in $A$
of an element in
$A_{d}$. In particular, if $\alpha_{1}$
and $\alpha_{2}$ are elements of $\widehat \mu (k)$ which are images of
elements
$\tilde \alpha_{1}$ and $\tilde \alpha_{2}$ of $\widehat \mu_{d} (k)$,
we denote by $J (\alpha_{1}, \alpha_{2})$ the image
of $J (\tilde \alpha_{1}, \tilde \alpha_{2})$ in $A$, which is 
independent from the choice of $d$, by Remark \ref{prop-quot}. For $\alpha$ in
$\widehat \mu (k)$ one defines similarly
$[\alpha (-1)]$ in $A$.

We will now consider the ring $U$ obtained
from the ring 
$A$ by adding all the Gauss sums motives associated to
$\widehat \mu (k)$.  This construction is strongly reminiscent of
Anderson's contruction of ``ulterior motives''
\cite{Anderson}.
\setcounter{equation}{0}

We define $U$ as the free $A$-module  with basis
$G_{\alpha}$, 
$\alpha$ in 
$\widehat \mu (k)$.
We define an $A$-algebra structure
on $U$ by putting the following relations~:
\begin{gather}
G_{1} = - 1\label{2.5.1}\\
G_{\alpha} G_{\alpha^{-1}} =
[\alpha (-1)] \, \LL \qquad \text{for} \quad \alpha \not=1\label{2.5.2}\\
G_{\alpha_{1}} G_{\alpha_{2}} =
J (\alpha_{1}, \alpha_{2}) \, G_{\alpha_{1}\alpha_{2}}
\qquad \text{for} \quad \alpha_{1}, \alpha_{2},
\alpha_{1}\alpha_{2}\not= 1.
\label{2.5.3}
\end{gather}

\begin{prop}The algebra $U$ is associative and commutative.
\end{prop}

\begin{proof}The commutativity is clear and the
associativity follows directly from Proposition \ref{formulaire}.
\end{proof}

\section{Motivic integrals of multiplicative characters}
\subsection{}Let $X$ be a $k$-variety, \ie
a separated and reduced $k$-scheme
of finite type.
We will denote by $\cL (X)$ the scheme of 
germs of
arcs on $X$. It is a scheme over $k$ and for any field extension
$k \subset K$ there is a natural bijection
$$\cL (X) (K) \simeq {\rm Mor}_{k-{\rm schemes}} (\Spec K [[t]], X)
$$ 
between the set of $K$-rational points of $\cL (X)$
and the set of 
germs of arcs with coefficients in  $K$ on $X$.
We will call
$K$-rational points of $\cL (X)$, for $K$
a field extension of $k$, arcs on $X$, and $\varphi (0)$ will be called 
the origin of the arc $\varphi$.
More precisely the
scheme $\cL (X)$ is defined as the projective limit
$$
\cL (X) := \varprojlim \cL_{n} (X)
$$
in the category of $k$-schemes of the schemes
$\cL_{n}(X)$ representing the functor
$$R \mapsto {\rm Mor}_{k-{\rm schemes}} (\Spec R [t] / t^{n+1} R[t], X)$$
defined on the category of $k$-algebras. (The existence of $\cL_{n}(X)$
is well known (cf. \cite{Arcs}) and the projective limit exists since 
the transition morphisms are affine.)
We shall denote by $\pi_{n}$ the canonical morphism, corresponding to 
truncation of arcs,
$$
\pi_{n} : \cL (X) \longrightarrow \cL_{n} (X).
$$
The schemes $\cL (X)$ and $\cL_{n} (X)$ will always be
considered with their reduced structure.

\subsection{}In \cite{Arcs}, we defined the boolean algebra $\BB_{X}$
of semi-algebraic subsets of $\cL (X)$.
We will refer to \cite{Arcs} for the precise definition,
but we
recall
that if $A$
is  a
semi-algebraic subset of $\cL (X)$,
the image $\pi_{n} (A)$ is constructible
in $\cL_{n} (X)$, and that for $f : X \rightarrow Y$
a morphism of $k$-varieties, the image by $f$ of any
semi-algebraic subset of $\cL (X)$ is
a
semi-algebraic subset of $\cL (Y)$. Both statements are consequences
from results by 
Pas \cite{Pas}.

Let $A$ be a semi-algebraic subset of $\cL (X)$. We call
$A$ {\it weakly stable at level} $n \in \NN$ if $A$ is a union of fibers
of $\pi_{n} : \cL (X) \rightarrow \cL_{n} (X)$. We call $A$ {\it weakly
stable}
if it stable at some level $n$. Note that weakly stable
semi-algebraic subsets form a boolean algebra.

Let $X$, $Y$ and $F$ be algebraic varieties over $k$,
and let
$A$, {\it resp.} $B$, be a constructible subset of $X$,
{\it resp.} $Y$. We say that
a map
$\pi : A \rightarrow B$ is a
{\it piecewise morphism}, if there exists a finite partition of $B$
in subsets $S$ which are
locally closed
in $Y$ such that $\pi^{- 1} (S)$ is locally closed in $X$ and
such that the restriction of $\pi$
to $\pi^{- 1} (S)$ is a morphism of $k$-varieties.
We say that
a map
$\pi : A \rightarrow B$ is a
{\it piecewise trivial fibration with fiber}
$F$, if there exists a finite partition of $B$ in subsets $S$ which are
locally closed
in $Y$ such that $\pi^{- 1} (S)$ is locally closed in $X$ and
isomorphic, as
a variety over $k$, to $S \times F$, with $\pi$
corresponding under the isomorphism to the projection
$S \times F \rightarrow S$. We say that the map $\pi$ is
a
{\it piecewise trivial fibration over} some constructible subset $C$ of
$B$,
if the restriction of $\pi$ to $\pi^{- 1} (C)$ is a piecewise 
trivial fibration onto $C$.
Let $X$ be an algebraic variety over $k$ of pure
dimension
$m$ and let $A$ be a semi-algebraic subset of $\cL (X)$. We call $A$
{\it stable at level} $n \in \NN$, if $A$ is weakly
stable at level $n$ and $\pi_{i + 1} (\cL (X)) \rightarrow \pi_{i} (\cL
(X))$
is a piecewise trivial fibration over $\pi_{i} (A)$ with fiber
$\AA^{m}_{k}$ for all $i \geq n$.

We call $A$ {\it stable} if it stable at some level $n$. Note that the
family of stable semi-algebraic subsets of $\cL (X)$ is closed under
taking finite intersections and finite unions. If $X$ is smooth, then $A$ is stable at level $n$
if it is weakly stable at level $n$.

\subsection{}We fix an integer $d \geq 1$ and assume $k$
contains all $d$-roots of unity.
Let $f : X \rightarrow \GG_{m, k}$ be a morphism of $k$-varieties.
For any character $\alpha$ of order $d$ of
$\mu_{d} (k)$, one may define an element $[X, f^{\ast}\cL_{\alpha}]$ of
$A_{d}$ as follows.

The morphism $[d] : \GG_{m, k} \rightarrow \GG_{m, k}$
given by $x \mapsto x^{d}$ is a Galois covering with Galois group
$\mu_{d} (k)$. We consider the fiber product
\begin{equation*}\begin{CD}
\widetilde X_{f, d} @>>>
X\\@VVV  @VV{f}V\\
\GG_{m, k} @>{[d]}>> \GG_{m, k}.
\end{CD}
\end{equation*}
The scheme $\widetilde X_{f, d}$ is endowed with an action of $\mu_{d} (k)$,
so we can define
$$[X, f^{\ast}\cL_{\alpha}]:= \chi_{c} (\widetilde X_{f, d}, \alpha).$$
More generally, if $X$ is constructible in some
$k$-variety and if $f : X \rightarrow \GG_{m, k}$
is a piecewise morphism, one may define
$[X, f^{\ast}\cL_{\alpha}] = \sum_{S \in \cS} [S,
f^{\ast}_{\vert S}\cL_{\alpha}]$, by taking an appropriate partition
$\cS$
of  $X$ into locally closed subvarieties, using the additivity of
$\chi_{c}$.

The following statement follows directly from the additivity of
$\chi_{c}$.

\begin{lem}\label{3.2.1}Let $X$ and $Y$ be constructible in some
$k$-varieties, and let
$f : X \rightarrow \GG_{m, k}$ and
$g : Y \rightarrow X$
be piecewise morphisms. Assume that $g$ is piecewise trivial
fibration with fiber $F$.
For any character $\alpha$ of order $d$ of
$\mu_{d} (k)$, the following holds
$$
[Y, (f \circ g)^{\ast} \cL_{\alpha}] =
\chi_{c} (F) [X, f^{\ast} \cL_{\alpha}].  \hfill \qed
$$
\end{lem}

\subsection{}\label{3.4}Let $X$ be an algebraic
variety
over $k$ of pure dimension $m$
and
let $f : X \rightarrow \AA^{1}_{k}$ be a morphism
of $k$-varieties. By the very definition
of semi-algebraic subsets, the set
$$\{{\rm ord}_{t} f = i \} := \{\varphi \in \cL (X) \bigm \vert
{\rm ord}_{t} f \circ \varphi = i\}$$
is a semi-algebraic subset of $\cL (X)$, for any
integer $i \geq 0$. One defines similarly the semi-algebraic subset
$\{{\rm ord}_{t} f > i \}$.
We will denote by $\bar f_{i}$ the mapping
$\{{\rm ord}_{t} f = i \} \rightarrow \GG_{m, k}$
which to a point $\varphi$ in $\{{\rm ord}_{t} f = i \}$
associates the constant term of the series
$t^{- i} (f \circ \varphi)$.
(Sometimes we shall use the same notation
$\bar f_{i}$ to denote the natural extension
$\{{\rm ord}_{t} f \geq i \} \rightarrow \AA^{1}_{k}$.)
Now let $W$ be a stable semi-algebraic subsets of $\cL (X)$
which is contained in $\{{\rm ord}_{t} f = i \}$
for some $i$. Choose an integer $n \geq i$
such that $W$ is stable at level $n$. 
The mapping $\bar f_{i}$ factors to a piecewise morphism
$$\bar f_{i \vert \pi_{n} (W)} : \pi_{n} (W) \rightarrow \GG_{m, k}.$$
Let $\alpha$ be a character of order $d$ of $\mu_{d} (k)$,
which we also view as an element of $\widehat \mu (k)$.
By Lemma \ref{3.2.1},
the virtual motive
$
[\pi_{n} (W), \bar f_{i \vert \pi_{n} (W)}^{\ast} \cL_{\alpha}]
\, 
\LL^{- (n + 1) m}
$
is independent of $n$.
So we may define
$$
\int_{W}^{\sim} \alpha ({\rm ac} f) d \mu$$
as the image of
$$[\pi_{n} (W), \bar f_{i \vert \pi_{n} (W)}^{\ast} \cL_{\alpha}]
\, \,
\LL^{- (n + 1) m}
$$
in $A$.

\begin{prop}Let $X$ be an algebraic
variety
over $k$ of pure dimension $m$
and
let $f : X \rightarrow \AA^{1}_{k}$ be a morphism
of $k$-varieties. For any $\alpha$ in $\widehat \mu (k)$,
there exists a unique map $$W \longmapsto
\int_{W} \alpha ({\rm ac} f) d \mu$$
from
$\BB_{X}$ to
$\widehat A$ 
satisfying the following three properties.
\begin{enumerate}
\item[] (3.4.2) \, If $W \in
\BB_{X}$ is stable and contained  in $\{{\rm ord}_{t} f = i \}$
for some $i$,
then
$\int_{W} \alpha ({\rm ac} f) d \mu$
coincides with the image of $\int_{W}^{\sim} \alpha ({\rm ac} f) d \mu$
in $\widehat A$ .
\item[] (3.4.3) \, If $W \in
\BB$ is contained in $\cL (S)$ with $S$ a closed subvariety of $X$ with 
${\rm dim} \, S < {\rm dim} \, X$, then
$\int_{W} \alpha ({\rm ac} f) d \mu = 0$.
\item[] (3.4.4) \, Let $W_{i}$ be in $\BB_{X}$ for each $i$ in $\NN$.
Assume that the
$W_{i}$'s are mutually disjoint and that
$W := \bigcup_{i \in \NN} W_{i}$ is semi-algebraic. Then
the series
$$\sum_{i \in \NN} 
\int_{W_{i}} \alpha ({\rm ac} f) d \mu$$
is convergent  in $\widehat A$ and converges
to $\int_{W} \alpha ({\rm ac} f) d \mu$.
\end{enumerate}
Moreover we have
\begin{enumerate}
\item[] (3.4.5) \, If $W$ and $W'$ are in $\BB_{X}$, $W \subset W'$,
and if $\int_{W'} \alpha ({\rm ac} f) d \mu$ belongs to the closure
$F^{m}
\widehat A$
of $F^{m} A$
in $\widehat A$,
then $\int_{W} \alpha ({\rm ac} f) d \mu$ 
belongs to
$F^{m} \widehat A$.
\end{enumerate}
\end{prop}
\begin{proof}Completely similar to the proof of Proposition 3.2 of
\cite{Arcs}.
\end{proof}

In the paper \cite{Arcs}, we also defined simple functions
$W \rightarrow \ZZ$, for
$W$ in $\BB_{X}$. A typical example of a simple function
is the following.
Consider a coherent sheaf of ideals $\cI$
on $X$
and  denote by 
${\rm ord}_t \cI$
the function
${\rm ord}_t \cI : \cL (X) \rightarrow \NN \cup \{+\infty\}$ given by
$\varphi \mapsto \min_{g} {\rm ord}_t g ( \varphi)$,
where the minimum is taken over all $g$ in the stalk $\cI_{\pi_{0}
(\varphi)}$
of $\cI$ at $\pi_{0}
(\varphi)$. The function ${\rm ord}_t \cI$ is a simple function.

Hence, for $W$ in $\BB_{X}$ and $\lambda : W \rightarrow \ZZ \cup
\{+ \infty\}$ a simple
function,
we can define
$$
\int_{W} \alpha ({\rm ac} f)
\LL^{- \lambda} d \mu := \sum_{n \in \ZZ}  \int_{W \cap \lambda^{-1}
(n)}
\alpha ({\rm ac} f) d \mu
\, \LL^{- n}
$$
in $\widehat A$, whenever the right hand side
converges,
in which case we say that $\alpha ({\rm ac} f) \LL^{- \lambda}$
is integrable on $W$. If the function
$\lambda$ is bounded from below, then
$\alpha ({\rm ac} f) \LL^{- \lambda}$
is integrable on $W$, because of (3.4.5).

For $Y$ a $k$-variety, we denote by
$\Omega^{1}_{Y}$ the sheaf of differentials on $Y$ and by
$\Omega^{d}_{Y}$ the
$d$-th exterior power of 
$\Omega^{1}_{Y}$. If $Y$ is smooth and $\cF$ is a coherent sheaf
on $Y$ together with a natural morphism
$\iota : \cF \rightarrow \Omega^{d}_{Y}$, we denote by $\cI (\cF)$ the sheaf of ideals on $Y$
which is locally generated by functions $\iota (\omega) / dy$ with $\omega$ a
local section of $\cF$ and $dy$ a local volume form on $Y$, and by
${\rm ord}_t \cF$ the simple function
${\rm ord}_t \cI (\cF)$.

\setcounter{theorem}{5}

\begin{prop}\label{3.3.6}Let $h : Y \rightarrow X$
be a proper birational morphism, with $Y$ smooth, let
$W$ be a semi-algebraic
subset of $\cL (X)$ and let $\lambda : \cL (X) \rightarrow
\NN$ be a simple function.
Then
$$
\int_W \alpha ({\rm ac} f) \LL^{- \lambda} d \mu = 
\int_{h^{ - 1} (W)} \alpha ({\rm ac} (f \circ h))
\LL^{- \lambda \circ h - {\rm ord}_t h^{\ast} (\Omega_{X}^{d})} d \mu.
$$
\end{prop}

\begin{proof}Follows directly from Lemma 3.4 of
\cite{Arcs}.
\end{proof}

\section{Motivic exponential integrals and the main result}

\subsection{}We set $\widehat U := U \otimes_{A}
\widehat A$.
We will also consider
the subring 
$A_{{\rm loc}}$
of $\widehat A$ generated by
the image of $A$ in
$\widehat A$
and the series $(1 - \LL^{-n})^{-1}$, $n \in \NN \setminus \{0\}$.
We denote by
$U_{{\rm loc}}$ the tensor product
$U \otimes_{A}
A_{{\rm loc}}$, which is naturally a subring of $\widehat U$.

\subsection{}Let $X$ be an irreducible
algebraic variety over $k$
and let
$f : X \rightarrow \AA^{1}_{k}$ be a morphism of
$k$-varieties.
Let $D$ be the divisor defined by $f = 0$ in
$X$. By a {\em very good}
resolution of
$(X, D)$, we mean a couple
$(Y, h)$ with
$Y$ is a smooth and connected $k$-scheme of finite type, $h : Y \rightarrow
X$
a proper morphism, such that the restriction
$h : Y \setminus h^{-1} (D \cup {\rm Sing} X) \rightarrow
X \setminus (D \cup {\rm Sing} X)$ is an isomorphism, and such that
the ideal sheaf  $\cI (h^{\ast} (\Omega^{d}_{X}))$
is invertible (this
last condition is of course irrelevant when
$X$ is smooth), and such that the union of
$(h^{-1} (D))_{\rm red}$ and the support of the divisor associated to
$\cI (h^{\ast} (\Omega^{d}_{X}))$ has only normal crossings as a
subscheme
of $Y$.
Such resolutions always exist by Hironaka's Theorem.

Let $E_{i}$, $i \in J$, be the irreduci\-ble (smooth) components
of
$(h^{-1} (D \cup {\rm Sing} X))_{\rm red}$.
Let $W$ be a reduced subscheme of $f^{-1} (0)$.
We call a very good
resolution $(Y, h)$
of
$(X, D)$ a 
very good
resolution
of
$(X, D, W)$ if $(h^{-1} (W))_{\rm red}$
is the reunion of some 
$E_{i}$'s.

Let $d$ be positive integer $\geq 1$. We will say
$d$ is big with respect to $(f, g, W)$ if
$$
\int_{\pi_{0}^{-1}(W)\cap \{{\rm ord}_{t}
f = i \}} \alpha ({\rm ac} f) \, \LL^{-{\rm ord}_t g}d \mu
= 0
$$
whenever the order of
$\alpha$  does not divide $d$.

\medskip

The following Proposition will be proven together with
Theorem \ref{MT}.
\begin{prop}\label{4.2.1}
Let $X$ be an irreducible
algebraic variety over $k$,
let
$f : X \rightarrow \AA^{1}_{k}$ and $g: X \rightarrow \AA^{1}_{k}$
be morphisms of
$k$-varieties. Let $W$ be a reduced subscheme of $f^{-1} (0)$.
Let $d$ be positive integer $\geq 1$. Assume
there exists
a very good
resolution $(Y, h)$
of
$(X, fg = 0, W)$ such that $d$ is a multiple of
the multiplicities of the $E_{i}$'s in the divisor
of $f \circ h$ on $Y$.
Then $d$ is big with respect to $(f, g, W)$.
\end{prop}

\begin{definition}Let $X$ be an irreducible
algebraic variety over $k$,
let
$f : X \rightarrow \AA^{1}_{k}$ and $g : X \rightarrow \AA^{1}_{k}$
be morphisms of
$k$-varieties. Let $W$ be a reduced subscheme of $f^{-1} (0)$.
For integers $i \geq 0$, we set \setcounter{equation}{1}
\begin{multline}\label{def-exp}
\int_{\pi_{0}^{-1}(W)} \exp (t^{- (i + 1)} f)
\, \LL^{-{\rm ord}_t g}d \mu
:=
\int_{\pi_{0}^{-1}(W) \cap \{{\rm ord}_{t} f > i \}} 
\, \LL^{-{\rm ord}_t g}d \mu\\
+
\sum_{\alpha \in \widehat \mu (k)} {\frac{1}{\LL - 1}}
G_{\alpha^{-1}}
\int_{\pi_{0}^{-1}(W) \cap \{{\rm ord}_{t}
f = i \}} \alpha ({\rm ac} f) \, \LL^{-{\rm ord}_t g}d \mu
\end{multline}
in $\widehat U$. 
The sum is finite since, by Proposition \ref{4.2.1}, there exists an 
integer $d$ which is big with respect to $(f, g, W)$.
\end{definition}

\begin{remark}The integral (\ref{def-exp}) belongs to
$U_{{\rm loc}}$ because of Proposition 5.1
in \cite{Arcs}
and the definitions in \ref{3.4}.
\end{remark}

\medskip
If
$f : X \rightarrow \AA^{1}_{k}$
and
$f' : X' \rightarrow \AA^{1}_{k}$ are morphisms of varieties,
we denote by $f \oplus f' : X \times X' \rightarrow \AA^{1}_{k}$
the morphism
given by composition of the morphism $(f, f')$ with the
addition morphism $\oplus : \AA^{1}_{k}
\times \AA^{1}_{k} \rightarrow \AA^{1}_{k}$.

\medskip

We can now state the main result of the paper.
\begin{MT}\label{MT}
Let $X$ and $X'$ be irreducible
algebraic varieties over $k$,
let
$f : X \rightarrow \AA^{1}_{k}$,
$g : X \rightarrow \AA^{1}_{k}$,
$f' : X' \rightarrow \AA^{1}_{k}$ and
$g' : X' \rightarrow \AA^{1}_{k}$
be morphisms of
$k$-varieties. Let $W$ (resp. $W'$)
be a reduced subscheme of $f^{-1} (0)$ (resp. $f'{}^{-1} (0)$).
For every $i \geq 0$, \setcounter{equation}{3}
\begin{multline}\label{4.2.4}
\int_{\pi_{0}^{-1}(W \times W')} \exp (t^{- (i + 1)} (f\oplus f'))
\, \LL^{-{\rm ord}_t gg'}d \mu
=\\
\Biggl(\int_{\pi_{0}^{-1}(W)} \exp (t^{- (i + 1)} f)
\, \LL^{-{\rm ord}_t g}d \mu\Biggr) \cdot
\Biggl(\int_{\pi_{0}^{-1}(W')} \exp (t^{- (i + 1)} f')
\, \LL^{-{\rm ord}_t g'}d \mu\Biggr).\,
\end{multline}
\end{MT}
\setcounter{equation}{5}

\begin{remark}Of course, by relations
(\ref{2.5.1})-(\ref{2.5.3}), (\ref{4.2.4}) is equivalent to
the following relations in
$A_{{\rm loc}}$, where,
for notational convenience, we will write
$[\alpha; f; g]$ for the integrand
$\alpha ({\rm ac} f) \, \LL^{-{\rm ord}_t
g}d \mu$,
\begin{multline}\label{4.2.6}
\int_{\pi_{0}^{-1}(W \times W')\cap \{{\rm ord}_{t}
f \oplus f' = i \}} [\alpha; f \oplus f'; 
gg'] =
\frac{1}{\LL - 1}
\sum_{\substack{\alpha_{1} \not=1 \, \alpha_{2} \not=1\\
\alpha_{1}\alpha_{2} = \alpha}}
J(\alpha_{1}^{-1}, \alpha_{2}^{-1}) \cdot \\
\Biggl(\int_{\pi_{0}^{-1}(W)\cap \{{\rm ord}_{t}
f = i \}} [\alpha_{1}; f ; g] \Biggr)
\cdot 
\Biggl(\int_{\pi_{0}^{-1}(W')\cap \{{\rm ord}_{t} 
f' = i \}} [\alpha_{2}; f' ; g']\Biggr)\\
-
\frac{1}{\LL - 1}
\Biggl(\int_{\pi_{0}^{-1}(W)\cap \{{\rm ord}_{t}
f = i \}} [1; f ; g]\Biggr)
\cdot  \Biggl(\int_{\pi_{0}^{-1}(W')\cap \{{\rm ord}_{t} 
f' = i \}} [\alpha; f' ; g']\Biggr)
+ \cdots\\
+
\Biggl(\int_{\pi_{0}^{-1}(W)\cap \{{\rm ord}_{t}
f > i \}} [1; f ; g]\Biggr)
\cdot 
\Biggl(\int_{\pi_{0}^{-1}(W')\cap \{{\rm ord}_{t} 
f' = i \}} [\alpha; f' ; g']\Biggr)
+ \cdots ,
\end{multline}
for $\alpha \not= 1$
(here $\cdots$ means ``same term as before with $(W, f, g)$ and $(W',
f', g')$
interchanged'')
and
\begin{multline}\label{4.2.7}
\int_{\pi_{0}^{-1}(W \times W')\cap \{{\rm ord}_{t}
f \oplus f' >  i \}} [1; f \oplus f'; 
gg'] - \\
\frac{1}{\LL -1}
\int_{\pi_{0}^{-1}(W \times W')\cap \{{\rm ord}_{t}
f \oplus f' =  i \}} [1; f \oplus f'; 
gg']=
\\
\sum_{\alpha \not= 1}
\frac{[\alpha (-1)] \LL}{(\LL - 1)^{2}}
\Biggl(\int_{\pi_{0}^{-1}(W)\cap \{{\rm ord}_{t}
f = i \}} [\alpha; f ; g]\Biggr)
\cdot
\Biggl(\int_{\pi_{0}^{-1}(W')\cap \{{\rm ord}_{t}
f' = i \}} [\alpha^{-1}; f' ; g']\Biggr)
\\
+
\Biggl(\int_{\pi_{0}^{-1}(W)\cap \{{\rm ord}_{t}
f >  i \}} [1; f ; 
g] - \frac{1}{\LL -1}
\int_{\pi_{0}^{-1}(W)\cap \{{\rm ord}_{t}
f  =  i \}} [1; f; 
g]\Biggr) \cdot 
\Biggl(\cdots\Biggr).
\end{multline}
Again, the sums are finite by Proposition \ref{4.2.1}.
\end{remark}

\subsection{Proof of Proposition \ref{4.2.1} and
Theorem \ref{MT}}
Applying
Pro\-po\-sition \ref{3.3.6} to very good resolutions of
$(X, fg = 0, W)$ and $(X', f'g' = 0, W')$ and using additivity
of $\chi_{c}$,
one reduces to the case when $X$ is smooth of dimension
$m$, $x_{1}, \dots, x_{m}$
are regular functions on $X$ inducing an \'{e}tale map
$X \rightarrow \AA^{m}_{k}$,
$f = u \prod_{i = 1}^{r} x_{i}^{n_{i}}$,
with $n_{i} >0$,
$g = v \prod_{i = 1}^{m} x_{i}^{m_{i}}$,
with $m_{i} \geq 0$, $u$ and $v$ are units, $W$ is the union
of the hypersurfaces
$x_{i} = 0$, for $i$ in $I \subset \{1, \ldots, r\}$,
and similarly with $'$
for
$X'$, $f'$, $g'$ and $W'$. We call this situation the DNC case.

We denote by $\tilde F^{n}_{d}$
the open subvariety of  $F^{n}_{d}$ defined by
$x_{i} \not=0$, $i = 1, \dots n$. It is stable under the 
$\mu_{d} (k)^{n}$-action.

\begin{prop}\label{4.3.1}Assume $d$
is big with respect to
$(f, g, W)$ and $(f', g', W')$.
For any $\alpha$ in $\widehat \mu (k)$ of order dividing $d$
and any integer $i \geq 0$,
the following relation holds in $A$,
\setcounter{equation}{0}
\begin{multline}
\int_{\substack{\pi_{0}^{-1}(W \times W')\cap \{{\rm ord}_{t}
f \oplus f' = i \}\\
\{{\rm ord}_{t}
f  = i \}\cap
\{{\rm ord}_{t}
f' = i \}
} } [\alpha; f \oplus f'; 
gg'] =
\frac{1}{\LL - 1}
\sum_{\alpha_{1}\alpha_{2} = \alpha}
\chi_{c} (\tilde F^{2}_{d}, (\alpha_{1}^{-1}, \alpha_{2}^{-1}))\\
\cdot \Biggl(\int_{\pi_{0}^{-1}(W)\cap \{{\rm ord}_{t}
f = i \}} [\alpha_{1}; f ; g] \Biggr)
\cdot 
\Biggl(\int_{\pi_{0}^{-1}(W')\cap \{{\rm ord}_{t} 
f' = i \}} [\alpha_{2}; f' ; g']\Biggr).
\end{multline}
Here we keep
the same notation for
$\chi_{c} (\tilde F^{2}_{d}, (\alpha_{1}^{-1}, \alpha_{2}^{-1}))$
and its image in $A$, and $\alpha_{1}$ and $\alpha_{2}$
are assumed to be of order dividing $d$.
\end{prop}

\begin{proof}By Proposition \ref{3.3.6}, one reduces to the DNC case.
For $i = 0$ the result is clear,
the domain of integration being empty,
so we may assume $i \geq 1$.
Now remark that $\pi_{0}^{-1}(W)\cap \{{\rm ord}_{t}
f = i \}$ is the disjoint union of semi-algebraic sets
$$
W_{\gamma} := \Bigl\{\varphi \Bigm | {\rm ord}_{t} x_{j} (\varphi) =
\gamma_{j} \Bigr\},
$$
for $\gamma = (\gamma_{1}, \dots, \gamma_{r})$,
with $\sum_{j = 1}^{r} \gamma_{j} n_{j} = i$ and $\gamma_{j} > 0$
when $j \in I$.
On $W_{\gamma}$ we may consider the function
$\bar x_{j} :  W_{\gamma} \rightarrow \GG_{m, k}$
which to a point 
$\varphi$ associates 
the constant term of the series
$t^{- \gamma_{j}}x_{j} (\varphi)$, for $1 \leq j \leq r$.
Set $$K = \Bigl\{j \in \{1, \dots, r\} \Bigm | \gamma_{j} > 0\Bigr\}.$$
The assumptions made imply that $K$ is not empty.
Now remark that on 
$W_{\gamma}$, $$\bar f_{i} = u \prod_{j \in K} \bar x_{j}^{n_{j}},$$
with $u$ a unit which is a function of the $\bar x_{j}$'s for
$j \notin K$ only.
Hence there exists $\ell > 0$ such that
$$W_{\gamma}
=
\pi_{\ell}^{-1} (Z \times \GG_{m, k}^{K}),
$$ with $Z$ a smooth variety,
and such that the restriction of $\bar f_{i}$ to $W_{\gamma}$
is equal to $$(u \prod_{j \in K} t_{j}^{n_{j}}) \circ \pi_{\ell},$$
with $u$ a unit on $Z$ and $t_{j}$ the canonical coordinate on the
corresponding
$\GG_{m, k}$ factor.
Let $a$ be the $gcd$ of the $n_{i}$'s.
By using an appropriate torus isomorphism of $\GG_{m, k}^{K}$
and changing $Z$,
one deduces that
$$W_{\gamma}
=
\pi_{\ell}^{-1} (Z \times \GG_{m, k}),
$$
with $Z$ a smooth variety,
and that
the restriction of $\bar f_{i}$ to $W_{\gamma}$
is equal to $(u t^{a}) \circ \pi_{\ell}$,
with $u$ a unit on $Z$ and $t$ the canonical coordinate on the
$\GG_{m, k}$ factor.
Now the result is a direct consequence of the next Lemma and the
following
Proposition. 
\end{proof}

\begin{lem}\label{4.3.2}With the above notations and assumptions,
if the order of $\alpha_{1}$ does
not divide $a$, then
$$\int_{W_{\gamma}} [\alpha_{1}; f ; g] = 0.$$
\end{lem}

\begin{proof}It is a direct consequence of Lemma 1.4.4 of
\cite{Motivic Igusa functions} (compare with the proof of
Proposition 2.2.2 (2) of
\cite{Motivic Igusa functions}).
\end{proof}

Remark that Proposition \ref{4.2.1} follows now directly from Lemma
\ref{4.3.2}.

\setcounter{equation}{2}
\begin{prop}Let $u_{1} :
X_{1} \rightarrow \GG_{m, k}$ and $u_{2} :
X_{2} \rightarrow \GG_{m, k}$ be morphisms of
algebraic varieties over $k$. Let $a ,b$ and $d$ be in $\NN \setminus
\{0\}$.
Assume $a$ and $b$ divide $d$.
Denote by
$(X_{1} \times \GG_{m, k} \times
X_{2} \times \GG_{m, k})^{0}$
the complement in 
$X_{1} \times \GG_{m, k} \times
X_{2} \times \GG_{m, k}$ of the divisor
of $u_{1} t_{1}^{a} + u_{2} t_{2}^{b}$.
For any character $\alpha$ of $\mu_{d} (k)$,
\begin{multline}\label{4.3.3}
[(X_{1} \times \GG_{m, k} \times
X_{2} \times \GG_{m, k})^{0}, 
(u_{1} t_{1}^{a} + u_{2} t_{2}^{b})^{\ast} \cL_{\alpha}]
= \\ (\LL - 1)
\sum_{\substack{
\alpha_{1} \alpha_{2} = \alpha\\
\alpha_{1}^{a} = 1 \, \, \alpha_{2}^{b} = 1
}}
\chi_{c} (\tilde F^{2}_{d}, (\alpha_{1}^{-1}, \alpha_{2}^{-1}))
[X_{1}, u_{1}^{\ast} \cL_{\alpha_{1}}]
[X_{2}, u_{2}^{\ast} \cL_{\alpha_{2}}].
\end{multline}
\end{prop}

\begin{proof}By definition the left hand side LHS of
(\ref{4.3.3}) is equal
to
$$
\chi_{c} \Bigl(
u_{1}(x_{1}) t_{1}^{a} + u_{2}(x_{2}) t_{2}^{b}
= w ^{d} \Bigm \vert
(x_{1}, x_{2}, t_{1}, t_{2}, w)
\in X_{1} \times X_{2} \times \GG_{m, k}^{3} \, \, ,
\alpha
\Bigr).
$$
Here the action of $\mu_{d} (k)$ is the standard one on the last $\GG_{m, k}$
and trivial on the other factors.
Hence, by Theorem
\ref{assertion},
\begin{multline*}
\text{LHS}
=
\chi_{c} \Bigl(u_{1}(x_{1}) = v_{1}^{a}\, , \,
u_{2}(x_{2}) = v_{2}^{b}\, , \,
v_{1}^{a}t_{1}^{a} + v_{2}^{b}t_{2}^{b}
= w ^{d} \Bigm \vert\\
(x_{1}, x_{2}, t_{1}, t_{2}, v_{1}, v_{2}, w)
\in X_{1} \times X_{2} \times \GG_{m, k}^{5} \, \, ,
(1, 1, \alpha)
\Bigr).
\end{multline*}
Here the group action is the action of 
$\mu_{a} (k) \times \mu_{b} (k) \times
\mu_{d} (k)$ which is
componentwise on the last three $\GG_{m, k}$ factors and trivial on
the others.
Making the toric change of variable
$$
(t_{1}, t_{2}, v_{1}, v_{2}, w)
\longmapsto
(T_{1} = v_{1} t_{1} w^{-d /a}, T_{2}= v_{2} t_{2} w^{-d /b}, v_{1},
v_{2}, w),
$$ this may be rewritten as
\begin{multline*}
\text{LHS}
=
\chi_{c} \Bigl(u_{1}(x_{1}) = v_{1}^{a}\, , \,
u_{2}(x_{2}) = v_{2}^{b}\, , \,
T_{1}^{a} + T_{2}^{b}
= 1 \Bigm \vert\\
(x_{1}, x_{2}, T_{1}, T_{2}, v_{1}, v_{2}, w)
\in X_{1} \times X_{2} \times \GG_{m, k}^{5} \, \, ,
(1, 1, \alpha)
\Bigr).
\end{multline*}
By
Proposition \ref{2.3.1} (3) and Lemma \ref{Kummer}
one deduces
$$
\text{LHS}
= (\LL - 1)
\sum_{\substack{\alpha_{1}  \in \widehat \mu_{a} (k)
\\
\alpha_{2} \in \widehat \mu_{b} (k)}}
[X_{1}, u_{1}^{\ast} \cL_{\alpha_{1}}]
[X_{2}, u_{2}^{\ast} \cL_{\alpha_{2}}] \,
\chi_{c} (\tilde F^{2}_{d}, (\alpha_{1}^{-1}, \alpha_{2}^{-1}, \alpha)).
$$
Here the 
$\mu_{a} (k) \times \mu_{b} (k) \times\mu_{d} (k)$-action on
$F^{2}_{d}$ is given
by
$$
(\xi_{1}, \xi_{2}, \xi_{3})
:
(T_{1}, T_{2}) \longmapsto
(\xi_{1} \xi_{3}^{- d/ a} T_{1}, \xi_{2} \xi_{3}^{- d/ b}T_{2}).
$$
The result follows because, by Proposition \ref{2.3.1} (4),
$\chi_{c} (\tilde F^{2}_{d}, (\alpha_{1}^{-1}, \alpha_{2}^{-1},
\alpha))$ is equal to $\chi_{c} (\tilde F^{2}_{d}, (\alpha_{1}^{-1},
\alpha_{2}^{-1}))$ (with the standard $\mu_{a} (k) \times \mu_{b} (k)$-action)
if $\alpha_{1} \alpha_{2} = \alpha$
and to $0$ otherwise.
\end{proof}

\begin{prop}\label{4.3.4}Let $X$ and $X'$ be irreducible
algebraic varieties over $k$,
let
$f : X \rightarrow \AA^{1}_{k}$,
$g : X \rightarrow \AA^{1}_{k}$,
$f' : X' \rightarrow \AA^{1}_{k}$ and
$g' : X' \rightarrow \AA^{1}_{k}$
be morphisms of
$k$-varieties. Let $W$ (resp. $W'$)
be a reduced subscheme of $f^{-1} (0)$ (resp. $f'{}^{-1} (0)$).
Let $i \geq 0$ be an integer.
\begin{enumerate}
\item[(1)]For any $\alpha \not= 1$ in
$\widehat \mu (k)$ and any integer $j > i$,
$$
\int_{\substack{\pi_{0}^{-1}(W \times W')\cap \{{\rm ord}_{t}
f \oplus f' = j \}\\
\{{\rm ord}_{t}
f  = i \}\cap
\{{\rm ord}_{t}
f' = i \}
} } [\alpha; f \oplus f'; 
gg'] = 0.
$$
\item[(2)]For any integer $j > i$,
\begin{multline*}
\int_{\substack{\pi_{0}^{-1}(W \times W')\cap \{{\rm ord}_{t}
f \oplus f' = j + 1 \}\\
\{{\rm ord}_{t}
f  = i \}\cap
\{{\rm ord}_{t}
f' = i \}
} } [1; f \oplus f'; 
gg'] = \\
\LL^{-1} \, \, 
\int_{\substack{\pi_{0}^{-1}(W \times W')\cap \{{\rm ord}_{t}
f \oplus f' = j \}\\
\{{\rm ord}_{t}
f  = i \}\cap
\{{\rm ord}_{t}
f' = i \}
} } [1; f \oplus f'; 
gg'].
\end{multline*}
\end{enumerate}
\end{prop}

\begin{proof}As before
one reduces to the DNC case and one uses the fact that
$$\pi_{0}^{-1}(W)\cap \{{\rm ord}_{t}
f = i \}$$ is the disjoint union of semi-algebraic sets
$$W_{\gamma} := \{\varphi  \, | \, {\rm ord}_{t} x_{j} (\varphi) =
\gamma_{j} \},$$
for $\gamma = (\gamma_{1}, \dots, \gamma_{r})$,
with $\sum_{j = 1}^{r} \gamma_{j} n_{j} = i$. We consider again
the functions
$\bar x_{j} :  W_{\gamma} \rightarrow \GG_{m, k}$.
On $W_{\gamma}$ we may write 
$
\bar f_{i} = v \prod_{1 \leq j \leq r} \bar x_{j}^{n_{j}}
$
with $v$ a unit, and similarly for $f'$.
Let $\varphi$ be a point in $W_{\gamma}$.
We write
$$
x_{j} (\varphi (t)) =
t^{\gamma_{j}} \bar x_{j} (1 + \sum_{k \geq 1} a_{k, j} t^{k}),
$$
for $1 \leq j \leq r$,
and
$$
x_{j} (\varphi (t)) =
\sum_{k \geq 0} a_{k, j} t^{k},
$$
for $r < j \leq m$.
Similar notation is used for $\varphi'$ in $W'_{\gamma'}$.
We may assume that $\gamma_{1} \geq 1$.
For $\ell > i$, the coefficient of $t^{\ell}$ in $f (\varphi (t)) +
f' (\varphi' (t))$ is equal to
\begin{equation*}\label{equ}\tag{3}
\sum_{j = 1}^{r} n_{j} a_{\ell - i, j} \bar f_{i}
+ P
+ 
\sum_{j' = 1}^{r'} n'_{j'} a'_{\ell - i, j'} \bar f'_{i} + P',
\end{equation*}
where $P$ (resp. $P'$) is a
polynomial in the variables
$a_{k,j}$ and $\bar x_{j}$ (resp. 
$a'_{k,j'}$ and $\bar x'_{j'}$), with $k \leq \ell -i$,
having as coefficients regular functions in
$\pi_{0} (\varphi)$ (resp. $\pi_{0} (\varphi')$).
Moreover, since $\gamma_{1} \geq 1$,
the polynomial $P$ does not involve the variable
$a_{\ell-1,1}$ (but might contain 
the variable
$a_{\ell-1,2}$ when $\gamma_{2} = 0$).

Let
$\Gamma_{\ell}$ be the locus of
${\rm ord}_{t} (f \oplus f') \geq \ell$ in 
$\pi_{\ell} (W_{\gamma} \times W'_{\gamma'}) \subset
\cL_{\ell} (X \times X')$,
and let $\Gamma^{+}_{\ell}$ be
the locus of ${\rm ord}_{t} (f \oplus f') > \ell$ in 
$\Gamma_{\ell}$.
From (\ref{equ}), for $\ell$ replaced by $i + 1$, \dots, $\ell$,
it follows that $a_{\ell-i, 1}$ does not appear in the equations
defining the variety $\Gamma_{\ell}$, and that 
$\Gamma^{+}_{\ell}$ is the hypersurface of $\Gamma_{\ell}$ defined by
equating (\ref{equ}) to zero. Taking the function
(\ref{equ}) as a new coordinate on $\Gamma_{\ell}$, instead of $a_{\ell
- 1, 1}$, we see that
$\Gamma_{\ell} \simeq \Gamma^{+}_{\ell} \times \AA^{1}_{k}$,
with the mapping $\overline{(f \oplus f')_{\ell}} : \Gamma_{\ell}
\setminus \Gamma_{\ell}^{+} \rightarrow \GG_{m, k}$
(which is given by (\ref{equ})) corresponding to the projection
of $\Gamma^{+}_{\ell} \times \GG_{m, k}$ onto the last factor. Assertion
(2) follows directly, and
assertion (1) is now a consequence from
Proposition \ref{2.3.1}  and Lemma \ref{Kummer}. 
\end{proof}

We are now able to conclude the proof of Theorem \ref{MT}.
Remark first that
$$
\chi_{c} (\tilde F^{2}_{d}, (\alpha_{1}, \alpha_{2})) =
\begin{cases}
J (\alpha_{1}, \alpha_{2}) & \text{if  $\alpha_{1}
\not = 1$ and $\alpha_{2} \not = 1$}\\
-1 & \text{if  $\alpha_{1}
\not = 1$ and $\alpha_{2} = 1$  or $\alpha_{1}
= 1$ and $\alpha_{2} \not= 1$}\\
\LL - 2
& \text{if  $\alpha_{1} = 1$ and $\alpha_{2} = 1$.}
\end{cases}
$$
If $\alpha \not= 1$, relation
(\ref{4.2.6}) follows directly from
Proposition \ref{4.3.1} and Proposition \ref{4.3.4} (1).
Assume now $\alpha = 1$.
We set $$a_{i} := \int_{\pi_{0}^{-1}(W)\cap \{{\rm ord}_{t}
f = i \}} [1; f ; g] \qquad
\text{and} \qquad
A_{i} := \int_{\pi_{0}^{-1}(W)\cap \{{\rm ord}_{t}
f > i \}} [1; f ; g]$$ and define similarly
$a'_{i}$ and $A'_{i}$.
For $k \geq i$, we also set
$$a_{i, k} := 
\int_{\substack{\pi_{0}^{-1}(W \times W')\cap \{{\rm ord}_{t}
f \oplus f' = k \}\\
\{{\rm ord}_{t}
f  = i \}\cap
\{{\rm ord}_{t}
f' = i \}
} } [1; f \oplus f'; 
gg'].$$

Let us denote by $\text{RHS}$ the right hand side
of (\ref{4.2.7}). 

Since, by Proposition \ref{formulaire},
$\chi_{c} (\tilde F^{2}_{d}, (\alpha, \alpha^{-1}))$ 
is equal to $-[\alpha (-1)]$ if $\alpha \not= 1$
and to $\LL - 2$
if $\alpha = 1$, we deduce from
Proposition \ref{4.3.1} the relation
\begin{multline*}\text{RHS} =
-\frac{\LL}{\LL - 1}
a_{i, i}
+
a_{i} a'_{i} + A_{i} A'_{i}
- \frac{1}{\LL - 1}a_{i} A'_{i} - \frac{1}{\LL - 1}a'_{i} A_{i}
=\\
-\frac{\LL}{\LL - 1}
\Bigl(a_{i}a'_{i} - \sum_{i < \ell}a_{i, \ell}\Bigr)
+
a_{i} a'_{i} + A_{i} A'_{i}
- \frac{1}{\LL - 1}a_{i} A'_{i} -  \frac{1}{\LL - 1}a'_{i} A_{i}.
\end{multline*}
The left hand side
$\text{LHS}$ of (\ref{4.2.7}) is
equal to
$$
A_{i}A'_{i}
+ \sum_{k \leq i < \ell}a_{k, \ell}
-\frac{1}{\LL - 1}
\Bigl(a_{i} A'_{i} + a'_{i}A_{i} + a_{i} a'_{i}
+ \sum_{k < i} a_{k, i} - \sum_{i < \ell} a_{i, \ell}
\Bigr).
$$
Hence we obtain
$$
\text{LHS} -
\text{RHS}
=
\sum_{k < i < \ell}a_{k, \ell} -
\frac{1}{\LL - 1} \sum_{k < i} a_{k, i} =
\sum_{k < i < \ell}a_{k, \ell} -
\frac{\LL^{-1}}{1 - \LL^{-1}} \sum_{k < i} a_{k, i}.
$$
The result now follows, since one deduces from
Proposition \ref{4.3.4} (2) that, for fixed $k$ and $i$,
with $k < i$,
$$\sum_{i < \ell}a_{k, \ell} =
\frac{\LL^{-1}}{1 - \LL^{-1}} \, \, a_{k, i}.\qed$$

\section{Motivic Thom-Sebastiani Theorem}

\subsection{}Let $B$ be any of the rings 
$A_{{\rm loc}}$, $\widehat A$, 
$U_{{\rm loc}}$, $\widehat U$.
We consider the ring of Laurent polynomials
$B[T, T^{-1}]$ and its localisation $B[T, T^{-1}]_{\rm rat}$
obtained by inverting the multiplicative family generated by
the polynomials $1 - \LL^{a} T ^{b}$, $a$, $b$ in $\ZZ$, $b \not = 0$.
Remark in this definition we could restrict to $b > 0$ or to $b < 0$.
Hence, by expanding denominators into formal
series, there are canonical embeddings of rings
$$
{\rm exp}_{T} : B[T, T^{-1}]_{\rm rat}
\longhookrightarrow B [T^{-1}, T]]$$
and
$${\rm exp}_{T^{-1}} : B[T, T^{-1}]_{\rm rat}
\longhookrightarrow B [[T^{-1}, T].
$$
Here $ B [T^{-1}, T]]$ (resp. $B [[T^{-1}, T]$) denotes the ring of
series
$\sum_{i \in \ZZ} a_{i} T^{i}$ with $a_{i} = 0$ for $i \ll 0$
(resp. $i \gg 0$).
By taking the difference ${\rm exp}_{T} - {\rm exp}_{T^{-1}}$
of the two expansions
one obtains an embedding 
$$
\tau :
B[T, T^{-1}]_{\rm rat} / B[T, T^{-1}]
\longhookrightarrow B [[T^{-1}, T]],
$$
where $B [[T^{-1}, T]]$ is the group of formal Laurent series
with coefficients in $B$.

Let
$\varphi =  \sum_{i \in \ZZ} a_{i} T^{i}$ and $\psi
=
\sum_{i \in \ZZ} b_{i} T^{i}$
be series in
$B [[T^{-1}, T]]$. We define their Hadamard product 
as the series
$$\varphi \ast \psi := \sum_{i \in \ZZ} a_{i} b_{i} \, T^{i}.$$

\begin{prop}\label{5.1.1}Let
$\varphi$ and $\psi$
be series in
$B [[T^{-1}, T]]$.
If they belong 
to the image of
$\tau$,
then their  Hadamard product $\varphi \ast \psi$
is also in the image of
$\tau$. 
\end{prop}

\begin{proof}Let $P_{1}$ and $P_{2}$ be in
$B [T, T^{-1}]_{\rm rat}$. From the formula
$$
\frac{1}{(1 - \LL^{a}T^{d})(1 - \LL^{b}T^{d})}
=
(\LL^{a} - \LL^{b})^{-1} 
\Bigl(\frac{\LL^{a}}{1 - \LL^{a}T^{d}} -
\frac{\LL^{b}}{1 - \LL^{b}T^{d}}
\Bigr)
$$
it follows
that there exists $d$ in $\NN \setminus \{0\}$ such that,
modulo $B [T, T^{-1}]$, both $P_{1}$ and $P_{2}$ are $B$-linear
combinations
of elements $T^{r} (1 - \LL^{a}T^{d})^{-k}$,
with $r \in \NN$, $a \in \ZZ$, $k \in \NN \setminus \{0\}$.
Thus, 
modulo $B [T, T^{-1}]$, both $\exp_{T} (P_{1})$ and
$\exp_{T} (P_{2})$ are $B$-linear combinations of elements of the form
$$
\varrho := \sum_{n \in \NN} f (n) \,  \LL^{na} \, T^{nd + r}
$$
in $B [T, T^{-1}]_{\rm rat}$, with $a \in \ZZ$,
$r \in \NN$, $r < d$, and $f$ a polynomial with coefficients in $\QQ$
such that $f (\ZZ) \subset \ZZ$.
We claim that
$$
\tau (\varrho) = \sum_{n \in \ZZ} f (n) \,  \LL^{na} \, T^{nd + r}.
$$
Indeed, if $f (n) = \binom{k + n - 1}{k - 1}$,
then $\varrho = T^{r} (1 - \LL^{a} T^{d})^{-k}$,
and an explicit calculation, using the relation
$$
\binom{k - m - 1}{k - 1} =
(-1)^{k - 1} \binom{m - 1}{k - 1},
$$
proves the claim in this special case. Hence the claim holds
for any $f$, considering $f$ as a linear combination
of such special $f$'s.
The Hadamard product of elements of the form $\varrho$
has again the same form. Thus the claim implies that the Hadamard
product commutes with $\tau$ for elements of the form $\varrho$, which
implies the result by the previous considerations.
\end{proof}

Let us denote by $B [[T]]_{\rm rat}$
the intersection of $B [[T]]$
with the image of ${\rm exp}_{T}$. It follows from the above proposition
that $B [[T]]_{\rm rat}$ is stable by Hadamard product.
Let $\varphi = {\rm exp}_{T} (P)$ be in 
$B [[T]]_{\rm rat}$. We denote by 
$\lambda (\varphi)$ the constant term
in the expansion of ${\rm exp}_{T^{-1}} (P)$.

\begin{prop}\label{5.1.2}Let $\varphi$ and
$\psi$ be series in 
$T B [[T]]_{\rm rat}$. Then 
$$
\lambda (\varphi \ast \psi)
=
- \, 
\lambda (\varphi) \cdot
\lambda (\psi).
$$
\end{prop}

\begin{proof}Let us remark that 
$\lambda (\varphi)$ only depends upon the class of $\varphi$
modulo additive translation by $T B [T]$. Hence we may assume there exists
$P$ and $Q$ in 
$B[T, T^{-1}]_{\rm rat}$ such that
${\rm exp}_{T} (P) = \varphi$,
${\rm exp}_{T} (Q) = \psi$,
${\rm exp}_{T^{-1}} (P)$ and 
${\rm exp}_{T^{-1}} (Q)$ belong to $B [[T^{-1}]]$. By Proposition
\ref{5.1.1} there exists $R$
in 
$B[T, T^{-1}]_{\rm rat}$ such that ${\rm exp}_{T} (R) = \varphi \ast
\psi$
and
${\rm exp}_{T^{-1}} (R) = -  {\rm exp}_{T^{-1}} (P)  {\rm exp}_{T^{-1}}
(Q)$.
The result follows.
\end{proof}

\subsection{}Let $X$ be an irreducible algebraic variety over $k$ of
pure dimension $m$
and
let $f : X \rightarrow \AA^{1}_{k}$ be a morphism. Let $W$ be
a reduced subscheme of $f^{-1} (0)$. 
We set
$$
E_{W, f} (T) =
\sum_{i > 0}
\Bigl[\int_{\pi_{0}^{-1} (W)}
\exp (t^{- (i + 1)} f) d\mu
\Bigr] T^{i}
$$
in $U_{{\rm loc}} [[T]]$. For any $\alpha$ in $\widehat 
\mu (k)$,
we set
$$
Z_{W, f, \alpha} (T) =
\sum_{i > 0}
\Bigl[\int_{\pi_{0}^{-1} (W) \cap
{\rm ord}_{t} f = i} \alpha ({\rm ac} f) d\mu
\Bigr] T^{i},
$$
in $A_{{\rm loc}} [[T]]$.
When $X$ is smooth,
$Z_{W, f, \alpha} (T)$ 
is equal to the natural image
in $A_{{\rm loc}} [[T]]$ of
$\int_{W} (f^{s}, \alpha)$, with the notation of
\cite{Motivic Igusa functions},
setting $T = \LL^{-s}$. Hence it follows from
Theorem 2.2.1 of \cite{Motivic Igusa functions}
that $Z_{W, f, \alpha} (T)$ belongs to
$A_{{\rm loc}} [[T]]_{{\rm rat}}$. This still holds
when $X$ is no more smooth by resolution of singularities and
Proposition \ref{3.3.6} (adapting the proof of Theorem 2.2.1 of
\cite{Motivic Igusa functions} in a straightforward way).

\medskip

We set
$$
\cS_{\alpha, W, f}^{\psi} := \frac{\LL^{m}}{1 - \LL} \lambda
(Z_{W, f, \alpha} (T))
$$
in $\widehat A$. When $X$ is smooth and
$\alpha$ is of order $d$, we defined
in Definition 4.1.2 of
\cite{Motivic Igusa functions} an element $S_{\alpha, x}$
of $A_{d}$, well defined modulo
$(\LL - 1)$-torsion,
for $x$ a closed point in $f^{-1} (0)$, which
corresponds
to the $\alpha$-equivariant part of
the motivic Euler characteristic of nearby cycles. Its image in
$\widehat A$ is just what we call now
$\cS_{\alpha, \{x\}, f}^{\psi}$. Remark there is no 
$(\LL - 1)$-torsion in $\widehat A$ (see also Remark \ref{separation}).
By Theorem 2.2.1 and Lemma 4.1.1 of \cite{Motivic Igusa functions},
$\cS_{\alpha, W, f}^{\psi}$ belongs to the image of
$A$, even when
$X$ is no more smooth, by resolution of singularities and
Proposition \ref{3.3.6}.

\medskip

To deal with vanishing cycles we set
$
\cS_{\alpha, W, f}^{\phi} = \cS_{\alpha, W, f}^{\psi}
$
for $\alpha \not = 1$,
and 
$
\cS_{\alpha, W, f}^{\phi} = \cS_{\alpha, W, f}^{\psi} - \chi_{c} (W)
$
for $\alpha = 1$.

\begin{prop}\label{5.2.1}The series
$E_{W, f} (T)$  belongs to 
$U_{{\rm loc}} [[T]]_{{\rm rat}}$ and
$$
\lambda ( E_{W, f} (T)) =
- \LL^{-m}
\sum_{\alpha \in \widehat \mu (k)} G_{\alpha^{-1}} 
\cS_{\alpha, W, f}^{\phi}.
$$
\end{prop}

\begin{proof}By the very definitions,
$$
E_{W, f} (T) = \frac{1}{\LL - 1}
\sum_{\alpha \in \widehat \mu (k)} G_{\alpha^{-1}}Z_{W, f, \alpha}
(T)
+ P (T)
$$
with $$P (T) = \sum_{i \geq 0} \int_{\pi_{0}^{-1} (W) \cap
\{{\rm ord}_{t} f > i\}}d\mu \, T^{i}
-
\int_{\pi_{0}^{-1} (W) \cap
\{{\rm ord}_{t} f > 0\}}d\mu
.$$
One may write $P (T)$ as             
$$
P (T) = (T - 1)^{-1} \Bigl(Z_{W, f, 1} (T) -
\int_{\pi_{0}^{-1} (W) \cap
\{{\rm ord}_{t} f > 0\}}d\mu\Bigr)-
\int_{\pi_{0}^{-1} (W)}d\mu.
$$
Since it
follows from Theorem 2.2.1 of
\cite{Motivic Igusa functions} that
$Z_{W, f, 1} (T)$ belongs to
$A_{{\rm loc}} [[T]]_{{\rm rat}}$ and that
${\rm exp}_{T^{-1}} (Z_{W, f, 1} (T)) $
belongs to $ A_{{\rm loc}} [[T^{-1}]]$, we deduce that
$\lambda (P (T)) = - \LL^{-m} \chi_{c} (W)$.
\end{proof}

\begin{MTS}\label{MTS}
Let $X$ and $X'$ be irreduci\-ble
al\-ge\-bra\-ic varieties over $k$ of pure dimension $m$ and $m'$.
Let
$f : X \rightarrow \AA^{1}_{k}$ and
$f' : X' \rightarrow \AA^{1}_{k}$ 
be morphisms of
$k$-varieties. Let $W$ (resp. $W'$)
be a reduced subscheme of $f^{-1} (0)$ (resp. $f'{}^{-1} (0)$).
Then 
$$
\sum_{\alpha} G_{\alpha^{-1}} \cS_{\alpha, W \times W' , f \oplus f'}^{\phi}
=
\Bigl(
\sum_{\alpha} G_{\alpha^{-1}} 
\cS_{\alpha, W, f}^{\phi}
\Bigr) \cdot
\Bigl(
\sum_{\alpha} G_{\alpha^{-1}} 
\cS_{\alpha, W', f'}^{\phi}\Bigr)
.$$
\end{MTS}

\begin{proof}By Theorem \ref{MT},
$E_{W \times W' , f \oplus f'} = E_{W, f }
\ast E_{W' , f'}$. The series
$E_{W, f } $ and $E_{W' , f'} $ having no constant term,
the result follows from Proposition
\ref{5.1.2}
and Proposition \ref{5.2.1}.
\end{proof}

\section{Hodge realization and the Hodge spectrum}

\subsection{}In this section we will assume $k = \CC$.
For $d \geq 1$, there is an embedding of $\widehat \mu_{d} (\CC)$ in
$\QQ / \ZZ$  given by $\alpha \mapsto a$
with $\alpha (e^{2 \pi i / d}) = e^{2 \pi i a}$. This gives an
isomorphism
$\widehat \mu (\CC) \simeq \QQ / \ZZ$. We denote by
$\gamma$ the section $\QQ / \ZZ \rightarrow
[0, 1)$.

A $\CC$-Hodge structure of weight $n$ is just a finite
dimensional bigraded vector space
$V = \bigoplus_{p + q = n} V^{p, q}$, or, equivalently,
a finite
dimensional vector space $V$ with decreasing filtrations $F^{\cdot}$ and
${\overline F}^{\cdot}$ such that $V = F^{p} \oplus {\overline F}^{q}$
when $p + q = n + 1$.  One defines similarly
a {\it rational} $\CC$-Hodge structure of weight $n$,
by allowing $p $ and $q$ to belong to $\QQ$ but still requiring
$p + q \in \ZZ$.

We denote by
$K_{0} ({\rm MHS}_{\CC})$ the Grothendieck group of the abelian category
of $\CC$-Hodge
structures
(it is also the Grothendieck group of the abelian category of complex
mixed Hodge
structures) and by 
$K_{0} ({\rm RMHS}_{\CC})$ the Grothendieck group of the abelian
category
of rational $\CC$-Hodge
structures.

The Hodge realization functor induces a morphism
$H : A \rightarrow K_{0} ({\rm MHS}_{\CC})$
which is zero on the kernel of the morphism
$A \rightarrow \widehat A$ by
Remark \ref{separation}. Hence we shall also write $H (A) $ for $A$ in
the image of $A$ by this morphism.

This morphism may be extended
to a morphism
$H : U \rightarrow K_{0} ({\rm RMHS}_{\CC})$
as follows.
For $p$ and $q$ in $\QQ$ with
$p + q$ in $\ZZ$, we denote by $H^{p, q}$ the class of the rank 1
vector space with bigrading $(p, q)$.
We set $H (G_{1}) = -1$ and
$H (G_{\alpha}) = - H^{1 - \gamma (\alpha),
\gamma (\alpha)}$ for $\alpha \not=1$.
This is compatible with the relations
\ref{2.5.1}--\ref{2.5.3} since, by a standard calculation (see \cite{Shioda}),
$$H (J_{\alpha_{1}, \alpha_{2}}) =
- H^{1 - (\gamma (\alpha_{1}) + \gamma (\alpha_{2}) -
\gamma (\alpha_{1} +
\alpha_{2})),
\gamma (\alpha_{1}) + \gamma (\alpha_{2}) -
\gamma (\alpha_{1} +
\alpha_{2})},$$
when $\alpha_{1} \not= 1$, $\alpha_{2} \not= 1$
and $\alpha_{1} \alpha_{2} \not= 1$.
Similarly as before, since $H$ vanishes on the kernel of the morphism
$U \rightarrow \widehat U$, we extend it to the image of this morphism.

\subsection{}For $X$ a complex
algebraic variety, we denote by
${\rm MHM} (X)$ the abelian category of mixed modules on $X$ constructed
by M. Saito in \cite{Sa1}, \cite{Sa3}.
In the definition of mixed Hodge modules it is required that the
underlying
perverse sheaf is defined over $\QQ$. To allow  some more
flexibility
we will also use the category ${\rm MHM'} (X)$ of bifiltred
$\cD$-modules on $X$ which are direct factors of objects of
${\rm MHM} (X)$ as bifiltred
$\cD$-modules. We denote
by $D^{b} ({\rm MHM} (X))$ and $D^{b} ({\rm MHM'} (X))$
the corresponding derived
categories.

Let $f : X \rightarrow \AA^{1}_{\CC}$ be a morphism.
We denote by $\psi^{H}_{f}$ and $\phi^{H}_{f}$
the nearby and vanishing
cycle functors for mixed Hodge modules as defined in \cite{Sa3}
and $T_{s}$
the semi-simple part of the monodromy operator.
One should note that $\psi^{H}_{f}$ and $\phi^{H}_{f}$
on mixed Hodge modules correspond to $\psi_{f}[-1]$ and $\phi_{f}[-1]$
on the underlying perverse sheaves.
If $M$ is a mixed Hodge module on $X$ 
we denote by $\psi^{H}_{f, \alpha} M$ the object
of
${\rm MHM'} (X)$ which corresponds to
the eigenspace of $T_{s}$ for the
eigenvalue
$\exp (2 \pi i \gamma (\alpha))$.
These definitions extend to the Gro\-then\-dieck group of the abelian
category
${\rm MHM'} (X)$.

For any object $K$ of $D^{b} ({\rm MHM} (X))$,
we denote by $\chi_{c} (X, K)$ the class of
$Rp_{!} (K)$ in $K_{0} ({\rm MHS}_{\CC})$, where $p$ is the projection onto
${\rm Spec} \, \CC$. Clearly this definition may be extended to
$D^{b} ({\rm MHM'} (X))$.

\begin{theorem}\label{6.2.1}Let $X$ be a smooth and connected
complex
algebraic variety of
dimension $m$, $f : X \rightarrow
\AA^{1}_{\CC}$ be a morphism and let
$i_{W} : W \hookrightarrow f^{-1} (0)$ be a
reduced subscheme of $f^{-1} (0)$. 
The following equalities hold
$$
H (\cS_{\alpha, W, f}^{\psi}) = (- 1)^{m - 1} \chi_{c} (W, i_{W}^{\ast} \psi_{f,
\alpha}^{H} \CC_{X}^{H} [m])
$$
and
$$
H (
\cS_{\alpha, W, f}^{\phi}) = (- 1)^{m - 1} \chi_{c} (W, i_{W}^{\ast} \phi_{f,
\alpha}^{H} \CC_{X}^{H} [m]).
$$
\end{theorem}

\begin{proof}When $W$ is a point the first equality is Theorem 4.2.1 of 
\cite{Motivic Igusa functions}. The proof of the general case
is completely similar. The second equality
follows directly from the first one.
\end{proof}

\begin{definition}Let $X$ be a smooth and connected
complex algebraic variety of
dimension $m$, $f : X \rightarrow
\AA^{1}_{\CC}$ be a morphism and let
$i_{W} : W \hookrightarrow f^{-1} (0)$ be a
reduced subscheme of $f^{-1} (0)$. 
We set
$$\tilde \chi_{c} (W, i_{W}^{\ast} \psi_{f}^{H} \CC_{X}^{H} [m])
= 
\sum_{\alpha} H (G_{\alpha^{-1}}) \chi_{c} (W, i_{W}^{\ast} \psi_{f,
\alpha}^{H} \CC_{X}^{H} [m])
$$
and
$$\tilde \chi_{c} (W, i_{W}^{\ast} \phi_{f}^{H} \CC_{X}^{H} [m])
= 
\sum_{\alpha} H (G_{\alpha^{-1}}) \chi_{c} (W, i_{W}^{\ast} \phi_{f,
\alpha}^{H} \CC_{X}^{H} [m]).
$$
\end{definition}

We deduce from Theorem \ref{MTS} and Theorem \ref{6.2.1}
the following Corollary.

\begin{cor}
Let $X$ and $X'$ be smooth and connected complex
algebraic varieties of pure dimension $m$ and $m'$.
Let
$f : X \rightarrow \AA^{1}_{\CC}$ and
$f' : X' \rightarrow \AA^{1}_{\CC}$ 
be morphisms of
algebraic varieties. Let $W$ (resp. $W'$)
be a reduced subscheme of $f^{-1} (0)$ (resp. $f'{}^{-1} (0)$).
Then 
\begin{multline*}\tilde \chi_{c} (W \times W', i_{W \times W'}^{\ast} \phi_{f
\oplus f'}^{H} \CC_{X \times X'}^{H} [m + m'])
=\\
-
\tilde \chi_{c} (W, i_{W}^{\ast} \phi_{f}^{H} \CC_{X}^{H} [m]) \, \cdot \,
\tilde \chi_{c} (W', i_{W'}^{\ast} \phi_{f'}^{H} \CC_{X'}^{H} [m']).
\qed
\end{multline*}
\end{cor}

Let $X$ be a smooth complex algebraic variety of pure dimension
$m$, let $f : X \rightarrow
\AA^{1}_{\CC}$ be a morphism of algebraic varieties and let
$x$ be
a closed point of $f^{-1} (0)$.
Let us recall the
definition of the spectrum ${\rm Sp} (f, x)$
given in \cite{St1} and \cite{Sa5}
(which  
differs from
that of \cite{St2} by multiplication by $t$).
Let $H$ be a complex
mixed Hodge structure with an
automorphism
$T$ of order dividing $d$. One defines 
the Hodge spectrum of $(H, T)$ as
${\rm HSp} (H, T) = \sum_{\alpha \in \frac{1}{d} \ZZ}  n_{\alpha}
t^{\alpha} \in
\ZZ [t^{- \frac{1}{d}}, t^{\frac{1}{d}}]$, with
$n_{\alpha} = {\rm dim \, Gr}^{p}_{F} H_{\lambda}$,
for $\lambda = \exp (2 \pi i \alpha)$ and $p = [\alpha]$,
where
$H_{\lambda}$ is the eigenspace of $T$ with eigenvalue $\lambda$,
and
$F$ is the Hodge filtration.
This definition extends to the Grothendieck group of the
abelian category of
complex mixed Hodge structures with an
automorphism
$T$ of order dividing $d$.
Note that ${\rm HSp} (H (k), T) = t^{-k} {\rm HSp} (H, T)$, where
$(k)$ is the Tate twist. We denote by
$\iota$ the $\ZZ$-algebra
automorphism of  $\ZZ [t^{- \frac{1}{d}}, t^{\frac{1}{d}}]$
defined by $\iota (t^{\frac{1}{d}})
= t^{- \frac{1}{d}}$.
Now one defines ${\rm Sp} (f, x)$ as
$${\rm Sp} (f, x) := t^{m} \iota \Bigl(\sum_{j \in \ZZ} (-1)^{j}  {\rm HSp}
(H^{j}i_{x}^{\ast}
\phi^{H}_{f} \CC_{X}^{H} [m], T_{s})\Bigr).
$$

\begin{cor}\label{spectrum}
Let $X$ and $X'$ be smooth and connected complex
algebraic varieties of pure dimension $m$ and $m'$.
Let
$f : X \rightarrow \AA^{1}_{\CC}$ and
$f' : X' \rightarrow \AA^{1}_{\CC}$ 
be morphisms of
algebraic varieties. Let $x$ and $x'$ be closed points in
$f^{-1} (0)$ and $f'{}^{-1} (0)$.
Then 
$${\rm Sp} (f \oplus f', (x, x')) =
{\rm Sp} (f, x) \, \cdot \, {\rm Sp} (f', x'). \qed
$$
\end{cor}

Corollary \ref{spectrum} was first proved by A. Varchenko in
\cite{Varchenko} when $f$ and $f'$ have
isolated singularities
(see also \cite{Scherk-Steenbrink}). The general
case is due to M. Saito (unpublished, but see
\cite{St2}, \cite{Sa6}).

\bibliographystyle{amsplain}

\end{document}